\documentclass[11pt]{article}
\usepackage{amssymb,amsmath,latexsym}
\usepackage{epsfig}

\newtheorem{thm}{Theorem}[section]

\newtheorem{cor}[thm]{Corollary}
\newtheorem{lemma}[thm]{Lemma}
\newtheorem{prop}[thm]{Proposition}
\newtheorem{defn}[thm]{Definition}

\numberwithin{equation}{section}

\def\1{{\bf 1}}
\def\ce{\mathcal{E}}
\def\cf{\mathcal{F}}

\def\pf{{\medskip\noindent {\bf Proof. }}}
\def\qed{{\hfill $\Box$ \bigskip}}

\def\sA {{\cal A}}  
  \def\sF {{\cal F}}
  
\def\sJ {{\cal J}}  \def\sL {{\cal L}}

 \def\bE {{\mathbb E}}

\def\bP {{\mathbb P}}  \def\bR {{\mathbb R}}

\def\nn{\nonumber}

\def\wt{\widetilde}
\def\wh{\widehat}

\def\E{{\mathbb E}}
\def\P{{\mathbb P}}

\def\bea{\begin{align*}}
\def\eea{\end{align*}}
\def\bee{\begin{equation}}
\def\eee{\end{equation}}

\def\eps{\varepsilon}

\def\R{{\mathbb R}}
\def\br{{\mathbb R}}
\def\b{{\beta}}

\def\wh{\widehat}

\makeatletter
\@addtoreset{equation}{section}

\makeatother

\begin{document}
\allowdisplaybreaks
\bibliographystyle{plain}

\title{\Large \bf
Two-sided estimates for the transition densities of symmetric Markov processes dominated by stable-like processes in $C^{1,\eta}$ open sets}

\author{{\bf Kyung-Youn Kim}\thanks{ 
This work was supported by a National Research Foundation of Korea (NRF) grant (No.~2009-0083521) funded by the Korean government (MSIP).}
\quad and \quad {\bf Panki Kim}\thanks{This work was supported by a National Research Foundation of Korea (NRF) grant (2013004822) funded by the Korean government (MEST).}
}
\date{}

\maketitle

\begin{abstract}
In this paper, we study sharp Dirichlet heat kernel estimates for a large class of symmetric Markov processes in $C^{1,\eta}$ open sets. 
The  processes are symmetric pure jump Markov processes with jumping intensity
$\kappa(x,y) \psi_1 (|x-y|)^{-1}  |x-y|^{-d-\alpha}$,
where $\alpha \in (0,2)$.
Here, $\psi_1$ is an increasing function on $[ 0, \infty )$, with $\psi_1(r)=1$ on $0<r \le 1$ and $c_1e^{c_2r^{\beta}} \le \psi_1(r) \le c_3 e^{c_4r^{\beta}}$ on $r>1$ for $\beta \in [0,\infty]$, and $ \kappa( x, y)$ is a symmetric function confined between two positive constants, with $|\kappa(x,y)-\kappa(x,x)|\leq c_5|x-y|^{\rho}$ for $|x-y|<1$  and $\rho>\alpha/2$. 
We establish two-sided estimates for the transition densities of such processes in $C^{1,\eta}$ open sets when $\eta \in (\alpha/2, 1]$.
In particular, our result includes (relativistic) symmetric stable processes and finite-range stable processes in $C^{1,\eta}$ open sets when $\eta \in (\alpha/2, 1]$. 
\end{abstract}


\bigskip

\section{Introduction}
Discontinuous Markov processes form a large class of stochastic processes containing stable-like processes and relativistic stable-like processes.
Recently, discontinuous Markov processes have often been used to simulate physical and economic systems that cannot be modeled by Gaussian processes (see \cite{KSZ, JW, M, OS, Ro}). 
Because of such importance in both theory and practice, there has been intense interest in studying discontinuous Markov processes. 

Throughout this paper we assume that $\beta \in [0, \infty]$, $\alpha\in (0, 2)$, and $d \ge 1$.
Let $\R^d$ be the $d$-dimensional Euclidean space and $dx$ be the $d$-dimensional Lebesgue measure in $\R^d$.
For $x\in \R^d$ and $r>0$, let $B(x, r)$ denote the open ball centered at $x$ with radius $r$.
The Euclidean distance between $x$ and $y$ will be denoted by $|x-y|$.
For two nonnegative functions $f$ and $g$, the notation $f\asymp g$ means that there are positive constants $c_1$ and 
$c_2$ such that $c_1g(x)\leq f (x)\leq c_2 g(x)$ in the common domain of definition for $f$ and $g$.
We will use the symbol ``$:=$,'' which is read as ``is defined to be.''
For any Borel set $A\subset \bR^d$, we will use ${\rm diam}(A)$ to denote its diameter and $|A|$ to denote its Lebesgue measure.

The infinitesimal generator $\sL$ of a discontinuous Markov process $Y=(Y_t, \P_x)_{t\ge 0, x\in \R^d}$ is a symmetric integro-differential operator, and under some mild assumptions 
the distribution $\P_x(Y_t \in dy)$ is absolutely continuous, for every $x  \in \R^d$ and $t>0$, with respect to Lebesgue measure in $\R^d$.
We will use $p(t, x, y)$ to denote the transition density of $Y$ so that  $\P_x(Y_t\in A)=\int_A p(t, x, y)dy$.
For any open subset $D\subset \bR^d$, we denote by  $Y^D$ the subprocess of $Y$ killed upon leaving $D$, and we use $p_D(t, x, y)$ to denote the transition density of $Y^D$.

The transition density $p_D(t, x, y)$ describes the distribution of the process $Y^D$.
Conversely, from an analytic viewpoint, $p_D(t, x, y)$ is also called a Dirichlet heat kernel of the operator $\sL$ on $D$, because it is a fundamental solution of $\partial_t u=\sL$ and $u = 0$ on  $D^c $.
Thus, obtaining sharp two-sided estimates of $p_D(t, x, y)$ is a fundamental problem in both analysis and probability theory.
However, it is not easy to obtain two-sided estimates of $p_D(t, x, y)$, especially near the boundary.
For Dirichlet heat kernel estimates for killed diffusions, see \cite{D2,Deb,DS} for the upper bound and \cite{Zq3} for the lower bound on bounded $C^{1,1}$ connected open sets. 

A prototype of discontinuous Markov processes is a (rotationally) symmetric $\alpha$-stable L\'evy process where $\alpha \in (0,2)$.
The  infinitesimal generator of a symmetric $\alpha$-stable L\'evy process is a fractional Laplacian   
$\Delta^{\alpha/2}=-(-\Delta)^{\alpha/2}$ 
that is a nonlocal operator. 
Recall that $\Delta^{\alpha/2}$ can be defined as 
\begin{align}\label{e:fLap}
\Delta^{\alpha/2}u(x)=\sA(d,-\alpha)\lim_{\varepsilon \to 0}\int_{\{y\in \R^d:|y-x|>\varepsilon\}}(u(y)-u(x))\frac{dy}{|x-y|^{d+\alpha}},
\end{align}
where $\Gamma$ is the Gamma function and $\sA(d,-\alpha)=\alpha 2^{\alpha-1}\pi^{-d/2}\Gamma(\frac{d+\alpha}{2})\Gamma(1-\alpha/2)^{-1}$.
Thus, it is a pure jump process and has a L\'evy density $y \to \sA(d,-\alpha)|y|^{-d-\alpha}$.
Chen et al.~\cite{CKS} obtained the Dirichlet heat kernel estimates for the symmetric $\alpha$-stable process $X$ in $C^{1,1}$ open sets. 

Another example of discontinuous Markov processes is a relativistic $\alpha$-stable process $X^m$ with mass $m>0$, which is a L\'evy process with a characteristic function given by
\[
\E_x\left[e^{i\xi\cdot(X_t^m-X_0^m)}\right]
=\exp\left(-t\left(\big(|\xi|^{2}+m^{2/\alpha}\big)^{\alpha/2}-m\right)\right)\mbox{ for every } x, \xi\in \R^d.
\] 
The corresponding infinitesimal generator is  $m-(m^{2/\alpha}-\Delta)^{\alpha/2}$.
In particular, for $\alpha=1$ the operator $m - \sqrt{m^2-\Delta}$ is called the free Hamiltonian corresponding to the quantization of the kinetic energy for a relativistic particle of mass $m$ (e.g., see~\cite{CMS, LY}). 
The L\'evy density of $X^m$ is 
\[
J^{m}(y)= \sA(d,-\alpha)|y|^{-d-\alpha}\psi(m^{1/\alpha}|y|) \quad \text{where} \quad
\psi (r):=  \int_0^\infty s^{\frac{d+\alpha}{2}-1} e^{-\frac{s}{ 4} -\frac{r^2}{ s} } \, ds.
\]
$\psi$ is decreasing and is a smooth function of $r^2$ satisfying $\psi(r)\le 1$ and
$\psi (r) \asymp e^{-r}(1+r^{(d+\alpha-1)/2} )$ on $[1,\infty)$  (see \cite[pp. 276--277]{CS4} for details).
Thus, $J^{m}(y)$ is dominated by the L\'evy density of the symmetric $\alpha$-stable process.
The approach developed in \cite{CKS} provides a guideline for establishing sharp two-sided heat kernel estimates for
other discontinuous L\'evy processes in open subsets of $\R^d$.
For example, two-sided Dirichlet heat kernel estimates for $X^m$ are discussed in~\cite{CKS2}.
Very recently two-sided Dirichlet heat kernel estimates were extended to a large class of symmetric L\'evy processes in~\cite{CKS6, CKS9}. 
 
In this paper, motivated by \cite{CKK3, CKS, CKS2} we consider a large class of symmetric Markov processes (not necessarily L\'evy processes)  whose jumping kernels are dominated by the kernel of the fractional Laplacian.
We establish the two-sided estimates for Dirichlet heat kernels of the generators of such Markov processes in (possibly unbounded) $C^{1,\eta}$ open sets $D$. 
When $D$ is $\R^d$, such a problem has been discussed in~\cite{KaSz, Sz1, Sz2}.
Our result extends the main results in \cite{CKS, CKS2} and provides far more. 

Let us now describe our assumptions and fix the notation simultaneously.
Let $\psi_1$ be an increasing function on $[0, \infty )$ with $\psi_1(r)=1$ for $0 < r\leq 1$, and let there be constants $\gamma_1,  \gamma_2> 0$ and $\beta \in [0, \infty]$ so that
\begin{equation}\label{eqn:exp}
  L_1  e^{\gamma_1r^{\beta}} \leq \psi_1 (r)\leq L_2  e^{\gamma_2r^{\beta}} \qquad
 \hbox{ for every }~ 1<r<\infty,
\end{equation}
for some constants $L_1, L_2>0$.
We define    
\begin{equation}\label{e:J3}
 j(r) =   \frac{1}{r^{d+\alpha}\psi_1 (r)} \quad r>0.
\end{equation}
We assume that $\kappa(x,y)$ is a positive symmetric function with 
\begin{equation}\label{e:conkappa1}
 L_3^{-1}\le \kappa(x,y) \le L_3, \quad x,y \in \R^d,
\end{equation}
and 
\begin{equation}\label{e:conkappa2}
  |\kappa(x,y)-\kappa(x,x)|{\bf 1}_{\{  |x-y|<1  \}}\leq L_4|x-y|^\rho, \quad x,y \in \R^d,
\end{equation}
where 
$\rho >\alpha/2$ and  $L_3, L_4$ are positive constants.
Let $J$ be a symmetric measurable function  on $\br^d\times \br^d \setminus \{x=y\}$ such that  \begin{equation}\label{e:J2} 
  J(x, y) :=  \kappa(x,y)   j(|x-y|) =
  \begin{cases}
  \kappa(x,y)  |x-y|^{-d-\alpha} \psi_1 (|x-y|)^{-1}  & \text{ if  } \beta \in [0, \infty),\\
  \kappa(x,y)   |x-y|^{-d-\alpha}{\bf 1}_{\{|x-y| \le1\}} & \text{ if  } \beta = \infty.
  \end{cases}
\end{equation}

For $u\in L^2(\br^d, dx)$, we define $\ce (u, u):=2^{-1}\int_{\br^d\times \br^d} (u(x)-u(y))^2 J(x, y) dx dy$.
Let $C_c(\br^d)$ denote the space of continuous functions with compact support in $\br^d$ and equipped with uniform topology.
We define
\begin{equation}\label{eqn:DF}
 \mathcal{D}(\ce):=\{f\in C_c (\br^d): \ce(f,f) <\infty\}.
\end{equation}
By \cite[Proposition 2.2]{CK2}, $(\ce, \cf)$  is a regular Dirichlet form on $L^2(\br^d,dx)$,  where
$\ce_1 (u, u):= \ce (u, u) +  \int_{\br^d} u(x)^2 dx$ and $\cf:=\overline{\mathcal{D}(\ce)}^{\ce_1}$. 
Hence, there is a Hunt process $Y$ associated with this on $\br^d$ (see~\cite{FOT}). 
 
Note that, since $j$ is decreasing and $J(x, y)\asymp j(|x-y|)$, we have
\[
\int_{B(x, r)} J(z, y)dz \ge \int_{B(x, r)\cap \{|z-y|\le|x-y|\}} c_1 j(|z-y|)dz \ge c_2 r^d j(|x-y|)\ge c_3r^{d} J(x, y)
\]
for every $r\le |x-y|/2$.
Thus, the Hunt process $Y$ associated with $(\ce, \cf)$ belongs to a subclass of the processes considered in~\cite{CKK3}.
Therefore, $Y$ is conservative and it has a H\"older continuous transition density $p(t, x, y)$ on $(0, \infty )\times \br^d\times \br^d$ with respect to the Lebesgue measure. 
 
The function $J$ is called the jumping intensity kernel of $Y$, because it gives rise to a L\'evy system for $Y$ describing the jumps of the process $Y$.
For any $x\in \bR^d$, stopping time $S$ (with respect to the filtration of $Y$), and nonnegative measurable function $f$ on $\bR_+ \times \bR^d\times \bR^d$ with $f(s, y, y)=0$ for all $y\in\bR^d$ and $s\ge 0$ we have
\begin{equation}\label{e:levy}
\E_x \left[\sum_{s\le S} f(s,Y_{s-}, Y_s) \right] = \E_x \left[ \int_0^S \left(\int_{\bR^d} f(s,Y_s, y) J(Y_s,y) dy \right) ds \right]
\end{equation}
(e.g., see \cite[Appendix A]{CK2}).
 
We first consider the estimate for the transition density $p(t,x,y)$ of $Y$ in $\R^d$.
Hereinafter, for $a, b\in \bR$, we have $a\wedge b:=\min \{a, b\}$ and $a\vee b:=\max\{a, b\}$.
For each  $a, T, \gamma>0$, we define a function $h_{a, \gamma, T}(t, r)$ on $(t,r) \in (0, T] \times [0, \infty)$ as
\begin{equation}\label{eq:qd}
h_{a, \gamma, T}(t, r):
=\begin{cases}
t^{-d /\alpha}\wedge  tr^{-d-\alpha}e^{-\gamma r^\b}&\text{ if } \b\in[0,1],\\
t^{-d /\alpha}\wedge  tr^{-d-\alpha}&\text{ if } \b\in (1, \infty] \text{ with } r<1,\\
t \exp\left\{-a \left( r \, \left(\log \frac{ T r}{t}\right)^{\frac{\b-1}{\b}}\wedge  r^\b\right) \right\}\qquad &\hbox{ if } \b\in(1, \infty) \text{ with } r \ge 1,\\
\left(t/(T r)\right)^{ar} 
&\hbox{ if }  \b=\infty \text{ with } r \ge 1.
\end{cases}
\end{equation}

Even though in \cite[Theorem 1.2]{CK2} and \cite[Theorems 1.2 and 1.4]{CKK3} two-sided estimates for $p(t, x, y)$ are stated separately for the cases $0<t\le 1$ and $t>1$, the constant 1 does not play any special role.
Thus, by the same proof, two-sided estimates for $p(t, x, y)$ hold for the case $0<t\le T$ and can be stated in an obvious way.

\begin{thm}\label{T:1.1}
Suppose that $Y$ is the symmetric pure jump Hunt process with the jumping intensity kernel $J$ defined in~\eqref{e:J2}. Then, the process $Y$ has a continuous transition density function $p(t, x, y)$ on $(0, \infty)\times \br^d\times \br^d$.
For each positive constant $T$, there are positive constants $C_1$, $c_1$, and $c_2\ge 1$ which depend on $\alpha, \beta, d, L_3, \psi_1, T$ such that for every $t \in (0,T]$ the function $p(t, x, y)$ has the following estimates:
\begin{align}\label{eq:M-THM1}
c_2^{-1}h_{c_1, \gamma_2, T}(t, |x-y|)\,\le\,
p(t,x,y)\,\le\, 
c_2\,h_{C_1, \gamma_1, T}(t, |x-y|).
\end{align}
\end{thm}

Note that, unlike those in \cite[Theorem 1.2]{CKK3}, the exponents $\gamma_1$ and $\gamma_2$ in Theorem~\ref{T:1.1} are explicit.  
When $\beta \in [0,1]$, the upper bound in \eqref{eq:M-THM1} comes from \cite[Theorem 2, Proposition 1]{KaSz}. 
We omit the proof of the upper bound in \eqref{eq:M-THM1} for $\beta \in [1, \infty]$, since the proof is the same, as mentioned above.
However, in Section~\ref{s:plbd} we give a detailed proof of the lower bounds in~\eqref{eq:M-THM1}. 

The goal of this paper is to obtain the sharp two-sided Dirichlet heat kernel estimates for $Y$ on $C^{1, \eta}$ open sets for $\eta\in(\alpha/2,1]$.
For any open set $D$, we use $\tau_D$ to denote the first exit time from $D$ by the process $Y$, and we use $Y^D$ to denote the process obtained by killing the process $Y$ upon exiting $D$.
By the strong Markov property, it can easily be verified that
$p_D(t, x, y):=p(t, x, y)-\E_x[ p(t-\tau_D, Y_{\tau_D}, y); t>\tau_D]$
is the transition density of $Y^D$.
Using the continuity and estimate of $p$, it is routine to show that $p_D(t, x, y)$ is symmetric and continuous (e.g., see the proof of Theorem~2.4 in~\cite{CZ}). 

Recall that an open set $D$ in $\bR^d$ (when $d\ge 2$) is said to be $C^{1,\eta}$  with $\eta\in(0, 1]$ if there exist a localization radius $ R>0 $ and a constant $\Lambda>0$ such that for every $z\in\partial D$ there exist a $C^{1,\eta}$-function $\phi=\phi_z: \bR^{d-1}\to \bR$ satisfying $\phi(0)=0$, $\nabla\phi (0)=(0, \dots, 0)$, $\| \nabla\phi \|_\infty \leq \Lambda$, $| \nabla \phi (x)-\nabla \phi (w)| \leq \Lambda |x-w|^{\eta}$ and an orthonormal coordinate system $CS_z$ of  $z=(z_1, \cdots, z_{d-1}, z_d):=(\wt z, \, z_d)$ with origin at $z$ such that $B(z, R )\cap D= \{y=({\tilde y}, y_d) \in B(z, R) \mbox{ in } CS_z: y_d > \phi (\wt y) \}$. 
The pair $( R, \Lambda)$ will be called the $C^{1,\eta}$ characteristics of the open set $D$.
Note that a $C^{1,\eta}$ open set $D$ with characteristics $(R, \Lambda)$ can be unbounded and disconnected, and the distance between two distinct components of $D$ is at least $R$.
By a $C^{1,\eta}$ open set in $\bR$ we mean an open set that can be written as the union of disjoint intervals so that the minimum of the lengths of all these intervals is positive and the minimum of the distances between these intervals is positive.

When $\beta \in (1, \infty]$, we need to make an assumption for $D$ in order to obtain the lower bound of $p_D(t,x,y)$.
We say that {\it the path distance in each connected component of $D$ is comparable to the Euclidean distance with characteristic $\lambda_1$} if for every $x$ and $y$ in the same component of $D$ there is a rectifiable curve $l$ in $D$ which connects $x$ to $y$ such that the length of $l$ is less than or equal to  $\lambda_1|x-y|$.
Clearly, such a property holds for all bounded $C^{1,\eta}$ open sets, $C^{1,\eta}$ open sets with compact complements, and connected open sets above graphs of $C^{1,\eta}$ functions.

We are now ready to state the main result of this paper.
Recall that $C_1$ is the constant in Theorem~\ref{T:1.1}. 
Let $\delta_D(x)$ be a distance between $x$ and $D^c$, and let
\begin{align}\label{e:dax}
\Psi(t,x):=\left(1\wedge\frac{\delta_D(x)^{\alpha/2}}{\sqrt{t}}\right).
\end{align}
\begin{thm}\label{t:main}
Suppose that $Y$ is the symmetric pure jump Hunt process with the jumping intensity kernel $J$ defined in~\eqref{e:J2}. 
Suppose that  $\eta\in(\alpha/2,1]$, $T>0$, and $D$ is a $C^{1,\eta}$ open set in $\bR^d$ with characteristics $( R, \Lambda)$.
Then, the transition density $p_D(t,x,y)$ of $Y^D$ has the following estimates.
\begin{description}
\item{\rm (1)}
There is a positive constant $c_1=c_1(\alpha, \b,  R, \Lambda, T, d , \psi_1, L_3,  L_4, \eta)$  such that for all $(t, x, y)\in (0, T]\times D\times D$ we have 
\begin{align*}
 p_D(t, x, y)
\leq \,c_1 \, \Psi(t,x) \Psi(t,y)
\begin{cases}
h_{ C_1\wedge\gamma_1, \gamma_1, T}(t, |x-y|/6)&\mbox{ if } \b\in[0,\infty),\\
h_{ C_1,  \gamma_1, T}(t, |x-y|/6)&\mbox{ if } \b=\infty.
\end{cases}
\end{align*}
\item{\rm (2)}
There is a positive constant $c_2=c_2(\alpha, \b,  R, \Lambda, T, d , \psi_1, L_3, L_4, \eta)$ such that for all $t \in (0, T]$ we have
\begin{align*}
p_D(t, x, y)\ge \,c_2 \,&\Psi(t,x)\Psi(t,y) 
\begin{cases}
t^{-d/\alpha}\wedge t|x-y|^{-d-\alpha} e^{-\gamma_2 |x-y|^\b}& \hskip -.1in\hbox{if } \b\in[0,1],\\
t^{-d/\alpha}\wedge t|x-y|^{-d-\alpha}& \hskip -.2in
\begin{array}{c}
\hbox{if  $\b\in(1, \infty)$ and $|x-y|<1$},\\
\hbox{ or  $\b=\infty$ and $|x-y| \le 4/5$}.
\end{array}
\end{cases}
\end{align*}
\item{\rm (3)}
Suppose in addition that the path distance in each connected component of $D$ is comparable to the Euclidean distance with characteristic $\lambda_1$.
Then, there are positive constants $c_i=c_i(\alpha, \b,  R, \Lambda, T, d, \psi_1, L_3, L_4, \eta, \lambda_1)$, $i=3,4$, such that if $x, y$ are in the same component of $D$ and $t \in (0, T]$, we have
\begin{align*}
p_D(t, x, y) \ge \,
c_3 \,\Psi(t,x) \Psi(t,y)  
\begin{cases}
h_{ c_4,  \gamma_2, T}(t, |x-y|)&\hbox{if $\b\in (1,\infty)$ and $|x-y|\ge 1$},\\
h_{ c_4,  \gamma_2, T}(t, 5|x-y|/4)&\hbox{if    $\b=\infty$ and $|x-y| \ge 4/5$}.
\end{cases}
\end{align*}
\item{\rm (4)}
If $\b\in (1, \infty)$, there is a positive constant $c_5=c_5(\alpha, \b, R, \Lambda, T, d, \psi_1,  L_3, L_4, \eta)$ such that for every $ x, y$ in the different components of $D$ with $|x-y|\ge 1$ and $t \in (0, T]$ we have  
\begin{align*}
  p_D(t, x, y) \ge\, 
c_5\,   \Psi(t,x) \Psi(t,y)\frac{t}{ |x-y|^{d+\alpha}}e^{-\gamma_2 (5|x-y|/4)^\b}.
\end{align*}
\item{\rm (5)}
Suppose in addition that $D$ is bounded and connected. 
Then, there are positive constants $c_i=c_i(\alpha, \b, R, \Lambda, T, d, \psi_1,  L_3, L_4, \eta, \text{\rm diam}(D)) $, $i=6,7$, such that for all  $(t, x, y)\in [T, \infty)\times D\times D$ we have
\[
c_6\, e^{- t \, \lambda^{D} }\, \delta_D (x)^{\alpha/2}\, \delta_D (y)^{\alpha/2}\,\leq\,
p_D(t, x, y) \,\leq\,
c_7\,e^{- t\, \lambda^{D}}\, \delta_D (x)^{\alpha/2} \,\delta_D (y)^{\alpha/2},
\]
where $-\lambda^D<0$ is the largest eigenvalue of the generator of $Y^D$.
\end{description}
\end{thm}

The cutoff value $5/4$ is not essential in the case $\beta=\infty$.
Further analysis reveals that for any $\eps>0$ we can choose $1+\eps$ as the cutoff value.
However, it seems that we cannot choose 1 as the cutoff value. 

If $D$ is a $C^{1,\eta}$ connected open set and the path distance in $D$ is comparable to the Euclidean distance, then by Theorem~\ref{t:main}(1)--(4) we can rewrite the two-sided estimates for $p_D(t, x, y)$.
\begin{cor}\label{c:dhke}
Suppose that $Y$ is the symmetric pure jump Hunt process with the jumping intensity kernel $J$ defined in~\eqref{e:J2}. 
Suppose further that $D$ is a $C^{1,\eta}$ connected open set with $\eta \in (\alpha/2,1]$ and that the path distance in $D$ is comparable to the Euclidean distance with characteristic $\lambda_1$. 
Then, for each  $T>0$ there exist $c_i=c_i(\alpha, \b, R, \Lambda, T, d, \psi_1,  L_3, L_4, \eta, \lambda_1)>0$, $i=1,2$, such that for every $t\in (0, T]$ we have 
\[
c_1^{-1} \, \Psi(t,x) \Psi(t,y) h_{ c_2, \gamma_2,  T}(t, 5|x-y|/4) \le  p_D(t, x, y)
\leq \,c_1 \, \Psi(t,x) \Psi(t,y) h_{ C_1\wedge\gamma_1, \gamma_1, T}(t, |x-y|/6).
\]
\end{cor}

The boundary Harnack principle for classical harmonic functions (for Brownian motion) describes how harmonic functions decay near the boundary of $D$.
This principle is important to studies of not only boundary value problems for partial differential equations but also the potential theory of Markov processes. 
The boundary Harnack principle has recently been generalized to a large class of discontinuous processes (see \cite{B, BKK, BKuK, G, KM2, KSV1, KSV7, SW}).

Unfortunately, the boundary Harnack principle does not hold for our process $Y$ when $\beta>1$ (see \cite{BKuK, KM2} for counterexamples).
This is one of the main difficulties in obtaining the boundary decay rate of $p_D(t,x,y)$.
In this paper, by using Dynkin's formula and the test function method, the key estimates for exit distributions are obtained directly.

Note that when $D$ is bounded, Theorem~\ref{t:main} gives the sharp estimates for $p_D(t, x, y)$ for all $t>0$, and the estimate for $p_D(t, x, y)$ has the same form as that obtained for symmetric stable processes in~\cite{CKS}.
Thus, by integrating the two-sided heat kernel estimates in Theorem~\ref{t:main} with respect to $t$ and following the proof of \cite[Corollary 1.2]{CKS}, the estimates for the Green function $G_D(x, y)=\int_0^\infty p_D(t, x, y)dt$ in \cite{CKS} can be extended to $C^{1,\eta}$ open sets. 
Since the proof is the same, we omit the proof.

\begin{cor}\label{C:1.2} 
Suppose that $Y$ is the symmetric pure jump Hunt process with the jumping intensity kernel $J$ defined in~\eqref{e:J2}. 
Suppose further that $\eta\in(\alpha/2,1]$ and $D$ is a bounded $C^{1,\eta}$ open set in $\bR^d$.
When $\beta=\infty$, we assume that $D$ is roughly connected. 
Then, on $D\times D$ we have
\[
G_D(x, y)\, \asymp\, \begin{cases} \displaystyle \frac{1} {|x-y|^{d-\alpha}}
\left(1\wedge \frac{  \delta_D(x)^{\alpha/2} \delta_D(y)^{\alpha/2}}{ |x-y|^{\alpha}}
\right)  \qquad &\hbox{when } d>\alpha ,  \\
\displaystyle \log \left( 1+ \frac{  \delta_D(x)^{\alpha/2} \delta_D
(y)^{\alpha/2}}{ |x-y|^{\alpha}
}\right)  &\hbox{when } d=1=\alpha , \\
\displaystyle\big( \delta_D(x)  \delta_D (y)\big)^{(\alpha-1)/2} \wedge \frac{
\delta_D(x)^{\alpha/2} \delta_D (y)^{\alpha/2}}{ |x-y|} &\hbox{when
} d=1<\alpha .
\end{cases}
\]
\end{cor}

The rest of this paper is organized as follows.
In Section~2, we first solve the Martingale-type problem for $Y$, which yields the Dynkin-type formula~\eqref{e:334}.
Then, in Theorem~\ref{L:2}, we give the key estimate for exit distributions. 
In Sections~3 and~5, we prove the lower bound estimates for $p_D(t,x,y)$.
In Section~3, we first consider the case $\delta_D(x)\wedge \delta_D(y)\ge t^{1/\alpha}$; that is, $x$ and $y$ are kept away from the boundary of $D$.
The result and our estimates for the exit distributions are used in Section~5 to prove the lower bound for all $x, y \in D$.
Section~4 contains the proof of the upper bound.
When $|x-y|<c$, we use Meyer's construction. 
Then, by using Lemma~\ref{L:4.1} twice, we prove the upper bound of $p_D(t,x,y)$ without using the lower bound of $p(t,x,y)$.
This enables us to write the bound of $p_D(t,x,y)$ in a compact form. 

 Throughout the rest of this paper, the positive constants $C_1, C_*, L_1, L_2, L_3, L_4, \gamma_1, \gamma_2$ can be regarded as fixed. 
In the statements of results and the proofs, the constants $c_i=c_i(a,b,c,\ldots)$, $i=1,2,3, \dots$, denote generic constants depending on $a, b, c, \ldots$, whose exact values are unimportant.
These are given anew in each statement and each proof.
The dependence of the constants on the dimension $d \ge 1$, on $\alpha \in (0,2)$, and on the positive constants $ L_1, L_2, L_3, L_4, \gamma_1, \gamma_2$ will not be mentioned explicitly.

\section{Estimates for exit distributions}

In this section we give some key estimates for exit distributions. 
First, we introduce an inequality that is used several times in this paper. 
\begin{lemma}\label{H:1}
Suppose that $\b\in[0,\infty)$.
For any $r_0>0$, there exists a positive constant $c=c( \b, r_0)$ such that
\begin{equation}\label{eH:1}
j(r)\le c j(2r)\qquad \hbox{for every}\quad  r\in(0, r_0].
\end{equation}
Moreover, \eqref{eH:1} holds for $\b=\infty$ if $r_0 \le 1/4$.
\end{lemma}
\pf
The result follows immediately from $L_2^{-1} e^{-\gamma_2 r^{\b}}r^{-d-\alpha}\le j(r)\le L_1^{-1} e^{-\gamma_1 r^{\b}}r^{-d-\alpha}$. \qed

For $\varepsilon \in (0,1/2)$, we define the operators $\sA_\eps$ and $\sA$ by 
\begin{align*}
\sA_\eps g(x):=\int_{\{y  \in \R^d:|y-x|>\varepsilon\}}(g(y)-g(x))J(y,x)dy \quad 
\quad \text{and}\quad \sA g(x):=  \lim_{\varepsilon \downarrow 0}\sA_\eps g(x)
\end{align*}
whenever these exist pointwise.
We use $C^{2}_c(\R^d)$ to denote the space of twice differentiable functions with compact support.
For every $g \in C_c^{2}(\R^d)$ and $r \in (\eps/2, 1]$ we have
\begin{align}
&\sA_\eps g(x)=\left(\int_{\{y  \in \R^d:r>|y-x|>\varepsilon\}} + \int_{\{y  \in \R^d:|y-x| \ge r\}}\right)(g(y)-g(x))\kappa (x,y) j(|x-y|)dy\nn\\
=&\kappa (x,x) \int_{\{y  \in \R^d:r>|y-x|>\varepsilon\}}(g(y)-g(x)) j(|x-y|)dy\nn\\
&+\int_{\{y  \in \R^d:r>|y-x|>\varepsilon\}}(g(y)-g(x))(\kappa (x,y)- \kappa (x,x))j(|x-y|)dy\nn\\
&+\int_{\{y  \in \R^d:|y-x|\ge r\}}(g(y)-g(x))\kappa (x,y) j(|x-y|)dy\nn\\
=&\kappa (x,x) \int_{\{y  \in \R^d:r>|y-x|>\varepsilon\}}\left(g(y)-g(x)-(y-x)
\cdot\nabla g(x)\right)j(|y-x|) dy\nn\\
&+\int_{\R^d} (g(y)-g(x))j(|x-y|)\left({\bf 1}_{\{r>|x-y|>\eps \}}(y)(\kappa (x,y)- \kappa (x,x))
+{\bf 1}_{\{|x-y| \ge r\}}(y) \kappa (x,y)\right)
dy\label{e:nag0}.
\end{align}
By \eqref{e:conkappa1} and \eqref{e:conkappa2} we have
\begin{align*}
&\left|(g(y)-g(x))j(|x-y|)\left({\bf 1}_{\{r>|x-y|>\eps \}}(y)(\kappa (x,y)- \kappa (x,x))
+{\bf 1}_{\{|x-y| \ge r\}}(y) \kappa (x,y)\right)\right|\\
\le &L_4{\bf 1}_{\{r>|x-y|>\eps \}}(y)|g(y)-g(x)| |x-y|^{-d-\alpha+\rho}+2L_3\|g\|_\infty{\bf 1}_{\{|x-y|>r\}}(y)  |x-y|^{-d-\alpha}
\end{align*}
and  $\rho > \alpha/2$. 
Thus, we see that $\sA g$ is well defined in $\R^d$ and that $\sA_\eps g$ converges to $\sA g$ locally uniformly in ${\R^d}$ as $\varepsilon\to 0$.
Furthermore, for every $r \in (0,1]$ we have
\begin{align}\label{e:nag}
&\sA g(x)=\kappa (x,x) \int_{\{y  \in \R^d:r>|y-x|\}}\left(g(y)-g(x)-(y-x)
\cdot\nabla g(x)\right)j(|y-x|) dy\nn\\
&+\int_{\R^d} (g(y)-g(x))j(|x-y|)\left({\bf 1}_{\{r>|x-y| \}}(y)(\kappa (x,y)- \kappa (x,x))
+{\bf 1}_{\{|x-y| \ge r\}}(y) \kappa (x,y)\right)
dy,
\end{align}
and
\begin{align*}
\|\sA g\|_\infty \le 
 c_1 \int_{\R^d} \left({\bf 1}_{\{1>|y| \}}(y) (|y|^{-d-\alpha+2} +
|y|^{-d-\alpha+\rho+1})
+{\bf 1}_{\{|y| \ge 1\}}(y) |y|^{-d-\alpha}\right)
dy < \infty.
\end{align*}

Next, we solve the Martingale-type problem for the operator $\sA$ on $C_c^{2}(\R^d)$ and show that the Dynkin-type formula in terms of $\sA$ is valid for every $f \in C_c^{2}(\R^d)$ (cf. \cite[Section 6]{GM}).
\begin{prop}\label{Dy}
For each $f \in C_c^{2}(\R^d)$ and $x\in \R^d$, there exists a $\P_x$-martingale $M^f_t$ with respect to the filtration of $Y$ such that $M^f_t=f(Y_t)-f(Y_0)-\int_0^{t}  {\sA} f(Y_s) ds$ is $\P_x$-a.s.
In particular,  for every $f \in C_c^{2}(\R^d)$ and 
any bounded open subset $U$ of $\R^d$ 
we have
\begin{equation} \label{e:334}
\E_x\int_0^{\tau_U}  {\sA} f(Y_t) dt =\E_x[f(Y_{\tau_U})]- f(x).
\end{equation}
\end{prop}

\pf 
We fix $f \in C_c^{2}(\R^d)$ and assume that the support of $f$ is a subset  of $B(0, R/2)$. 
We use a strict version of Fukushima's decomposition \cite[Theorem 5.2.5]{FOT}.
First, it is clear from \eqref{eqn:DF} that $ f \in  \sF$.
The energy measure $\mu_{\langle f \rangle}$ of $f$ has the density $\Gamma(f)(x)= \int_{\br^d} (f(x)-f(y))^2 J(x, y) dy$.
Note that $\|\Gamma(f)\|_\infty < \infty$ and that $|\Gamma(f)(x)|\le c_1 |x|^{-d-\alpha}$ for $x \in  B(0, R)^c$.
Thus, $\mu_{\langle f \rangle}(\R^d)<\infty$.

Now, by Fubini's theorem and the dominated convergence theorem, for any $g\in C^2_c(\R^d)$ we have
\begin{eqnarray*}
&&\ce(f,g)=\frac12\lim_{\varepsilon\downarrow0}\int_{\{(x, y)\in {\R^d}\times {\R^d},\
|y-x|>\varepsilon\}}(g(y)-g(x))(f(y)-f(x))J(y,x)\, dx\, dy\\
&&=
-\lim_{\varepsilon\downarrow0}\int_{\R^d} g(x)
\left(\int_{\{y \in \R^d:|y-x|>\varepsilon\}}(f(y)-f(x))J(y,x) dy \right) dx\,=\,
-\,\int_{\R^d} g(x)\sA f(x)\, dx.
\end{eqnarray*}

We recall from \cite{FOT} that $S_0$ is the collection of positive Radon measures of finite energy integrals and 
\[
S_{00}:=\{ \mu \in S_0 :\mu(\R^d) < \infty, \sup_{x \in \R^d}\int_{\R^d} \int_0^\infty e^{- t} p(t,x,y)dt \mu(dy) < \infty \} 
\]

Let $\nu:=\nu_+-\nu_-$, where $\nu_+(dx):=-{\bf 1}_{\{\sA f(x)<0\}}\sA f(x)dx$ and $\nu_-(dx):={\bf 1}_{\{\sA f(x) \ge 0\}}\sA f(x)dx$, so that $\ce(f,g)=\int_{\R^d} g(x) \nu(dx)$.
Note that $\|\sA f\|_\infty < \infty$ and that $|\sA f(x)|\le c_2 |x|^{-d-\alpha}$ for $x \in B(0, R)^c$.
Thus ,$|\nu|(\R^d) <\infty$.
Moreover, clearly 
\begin{eqnarray*}
\sup_{x \in \R^d} \int_{\R^d} \int_0^\infty e^{- t} p(t,x,y)dt |\nu|(dy)
\le \|\sA f\|_\infty \int_0^\infty e^{- t}dt  <\infty. 
\end{eqnarray*}
Thus, $\nu_+$ and $\nu_-$ are in $S_{00}$.
Since $-\int_0^{t} {\bf 1}_{\{\sA f(X_s)<0\}} {\sA} f(X_s) ds$  and $\int_0^{t} {\bf 1}_{\{\sA f(X_s) \ge 0\}} {\sA} f(X_s) ds$ are  positive continuous additive functionals in the strict sense with Revuz measures $\nu_{+}$ and $\nu_{-}$, respectively, 
upon applying \cite[Theorem 5.2.5]{FOT} we conclude that for every $x \in \R^d$ we have 
\begin{align*}
f(Y_t)-f(Y_0)&=M^{f}_t + \int_0^{t} {\bf 1}_{\{\sA f(x) \ge 0\}} {\sA} f(Y_s) ds+ \int_0^{t} {\bf 1}_{\{\sA f(x) < 0\}} {\sA} f(Y_s) ds\\
&=M^{f}_t + \int_0^{t}  {\sA} f(Y_s) ds,
\end{align*}
where $M^{f}_t$ is a $\P_x$-martingale additive functional in the strict sense with Revuz measure $\mu_{\langle f \rangle}$. \qed 

Using \eqref{e:334}, we prove the following lemma, which is used several times in Section~4. 
\begin{lemma}\label{L:2.3}
For every $a \in (0,1]$, 
there exists a positive constant $c=c(a)$ such that, for any $\beta \in [0, \infty]$, any $r \in (0,1]$, and any open sets $U$ and $D$ with $B(0, a r ) \cap D \subset U  \subset  D$, we have
\[
\P_x\left(Y_{\tau_U} \in D\right) \,\le\, c\,r^{-\alpha}\, \E_x[\tau_U],
\qquad x \in D\cap B(0, a r/2).
\]
\end{lemma}

\pf 
For fixed $a \in (0,1]$, 
we take a sequence of radial functions $(\phi_m)_{m \ge 1}$ in $C^{\infty}_c(\R^d)$ such that $0\le \phi_m\le 1$, with
\[
\phi_m(y)=\left\{
\begin{array}{lll}
0, & \text{if } |y|< a/2 \text{ or } |y|>m+2,\\
1, & \text{if } a\le |y|\le m+1,\\
\end{array}
\right.
\]
and
\[
\sup_{m\ge1}\left (\sum^d_{i=1}\|\frac{\partial}{\partial y_i}
\phi_m\|_\infty+ \sum^d_{i, j=1}\|\frac{\partial^2}{\partial y_i\partial y_j}
\phi_m\|_\infty \right) < c_1=c_1(a)< \infty.
\]
For any $r \in (0,1]$, define $\phi_{m, r}(y)=\phi_m(\frac{y}{r})$ so that
$0\le \phi_{m, r}\le 1$,
\begin{equation}\label{e:2.111}\phi_{m, r}(y)=
\begin{cases}
0, & \text{if }|y|<ar/2 \text{ or } |y|>r(m+2)\\
1, & \text{if } ar\le |y|\le r(m+1), 
\end{cases}
\end{equation}
and such that
\begin{equation}\label{e:2.11}
\sup_{m \ge 1 } \sum^d_{i=1}\left\|\frac{\partial}{\partial y_i} \phi_{m,
r}\right\|_\infty \,<\, c_1\, r^{-1}
\quad \text{and} \quad  \sup_{m \ge 1 } \sum^d_{i,
j=1}\left\|\frac{\partial^2}{\partial y_i\partial y_j} \phi_{m,
r}\right\|_\infty \,<\, c_1\, r^{-2}.
\end{equation}
Using \eqref{e:conkappa1}, \eqref{e:conkappa2}, \eqref{e:2.11}, and the assumption 
that $\rho >\alpha/2$, for every $x \in \R^d$, $r \in (0,1]$, and $m \ge 1$ we have
\begin{align}
 &
\left|\kappa (x,x) \int_{\{y  \in \R^d:r>|y-x|\}}\left(\phi_{m,r}(y)-\phi_{m,r}(x)-(y-x)
\cdot\nabla \phi_{m,r}(x)\right)j(|y-x|) dy
 \right|\nn\\
&+\int_{\R^d} |\phi_{m,r}(y)-\phi_{m,r}(x)|j(|x-y|)\left({\bf 1}_{\{r>|x-y| \}}(y)|\kappa (x,y)- \kappa (x,x)|
+{\bf 1}_{\{|x-y| \ge r\}}(y) \kappa (x,y)\right)
dy\nn\\
&\le \frac{c_2}{r^2}\int_{\{|x-y|< r \}}  |x-y|^{-d-\alpha+2} dy   + 
\frac{c_2}{r}\int_{\{|x-y|< r \}}  |x-y|^{-d-\alpha+1+\rho} dy   +c_2
\int_{\{|x-y| \ge r\}}|x-y|^{-d-\alpha} dy \nn\\
& \le {c_3}(r^{-\alpha}+ r^{-\alpha+\rho})
 \le {2c_3}r^{-\alpha}\,  \label{e:2.12}
\end{align}
for some $c_3=c_3(a)>0$.
Now, by combining \eqref{e:nag}, \eqref{e:334}, \eqref{e:2.111}, and \eqref{e:2.12}, we find that for any $x\in D\cap B(0, ar/2)$ we have
\begin{align*}
&\P_x\left(Y_{\tau_U} \in \{ y \in D: ar \le |y| <(m+1)r \}\right)
=\E_x\left[\phi_{m, r} \left(Y_{\tau_U}\right): Y_{\tau_U} \in \{ y
\in D: ar  \le |y| <(m+1)r
\}\right]\\
&\le \E_x\left[\phi_{m, r} \left(Y_{\tau_U}\right)\right]=
\E_x\left[ \int_0^{\tau_U}   
{\sA} \phi_{m, r}(Y_t)dt \right] \le
\|{\sA}\phi_{m, r}\|_\infty\,  \E_x[\tau_U]
\le 2c_{3} r^{-\alpha}\E_x[\tau_U].
\end{align*}
Therefore, since $B(0, a r ) \cap D \subset U$, we obtain
\begin{align*}
\P_x\left(Y_{\tau_U} \in D\right)= \lim_{m\to
\infty}\P_x\left(Y_{\tau_U} \in \{ y \in D: a r \le |y| <(m+1)r
\}\right) \,\le\,
2c_3\,r^{-\alpha}\E_x[\tau_U].
\end{align*}
\qed

For the remainder of this section we assume that $\eta\in(\alpha/2,1]$ and that $D$ is a $C^{1, \eta}$ open set with $C^{1, \eta}$ characteristics $(R, \Lambda)$.
Without loss of generality, we assume that $R\leq 1$ and $\Lambda \geq 1$. 
For each fixed $Q \in \partial D$ and for every $r \le R$ we define
\begin{align}\label{e:h}
h_{Q,r}(y) :=  \delta_{D}(y)^{\alpha/2} {\bf 1}_{D\cap B(Q, r)}(y).
\end{align}

We next establish two lemmas that are used to obtain the key estimates for exit distribution. 
The next lemma and its proof are similar to \cite[Lemma 2.3]{CKS3} and \cite[Lemma 3.7]{KM2} and their proofs.
We provide the proof here for completeness. 
Recall that $\Delta^{\alpha/2}$ is defined in~\eqref{e:fLap}.

\begin{lemma}\label{l:frac}
There exists a positive constant $c=c(\eta, R, \Lambda)$ independent of $Q \in \partial D$ such that $\Delta^{\alpha/2} h_{Q,R/2}$ is well defined in ${D}\cap B(Q, R/8)$ and
\[
| \Delta^{\alpha/2} h_{Q,R/2}(x)| \le  c \quad \text{ for all }
x \in {D}\cap B(Q, R/8)\, .
\]
\end{lemma}

\pf
Since the case of $d=1$ is easier, we give the proof only for $d\ge 2$. 

Let $h(\cdot):=h_{Q,R/2}( \cdot)$.
Fix $x\in D\cap B(Q, R/8)$ and let $z_x\in \partial D$ such that $\delta_D(x)=|x-z_x|$.
Let $\phi$ be a  $C^{1, \eta}$ function and $CS=CS_{z_x}$ be an orthonormal coordinate system with $z_x$ chosen so that $x=(\wt{0}, x_d)$, $B(0, R)\cap D=\{y=(\wt{y}, y_d) \mbox{ in } CS:y \in B(0, R), y_d>\phi(\wt{y})\}$, 
$\phi(\wt 0)=0$, $ \nabla\phi(\wt 0)=(0, \dots, 0)$, $\| \nabla \phi \|_\infty \leq \Lambda$, and $|\nabla \phi (\wt y)-\nabla \phi_Q (\wt z)| \leq \Lambda |\wt y-\wt z|^\eta$.
We fix the function $\phi$ and the coordinate system $CS$, and we define a function $h_x(y)=\delta_{H^+}(y)^{\alpha/2}$, where  $H^+=\{y=(\wt{y}, y_d)$ in $CS:y_d>0\}$ is the half space in $CS$. 

We define $\wh \phi :B(\wt{0}, R)\to \bR$ by $\wh \phi(\wt{y}):=2\Lambda|\wt{y}|^{\eta+1}$.
Since $\nabla\phi (0)=0$, by the mean value theorem we have 
$-\wh \phi(\wt{y})\le\phi(\wt{y})\le \wh \phi(\wt{y})$  for any $y\in D\cap B(x, R/8)$.
Since $\Delta^{\alpha/2} h_x(y)=0$ for any $y \in H^+$(see Lemma~2.1 of \cite{CKS3}), it is enough to show that $\Delta^{\alpha/2}(h- h_x)(x)$ is well defined and that there exists 
a constant $c_1=c_1(\eta, R, \Lambda)>0$ independent of  $x\in D\cap B(Q, R/8)$ and $Q \in \partial D$  such that
\begin{align}\label{e:claims}
\int_{D\cup H^+}\frac{|h(y)-h_x(y)|}{|y-x|^{d+\alpha}} dy\le c_1< \infty.
\end{align}

Let $A:=\{y\in (D\cup H^+)\cap B(x, R/8):-\wh \phi(\wt{y})\le y_d\le \wh \phi(\wt{y})\}$
and $E:=\{y\in B(x, R/8): y_d>\wh \phi(\wt{y})\}$. 
We prove \eqref{e:claims} by showing that $I+II+III\le c_1$, where
\begin{align*}
I:=\int_{B(x, R/8)^c} \frac{h(y)+h_x(y)}{|y-x|^{\alpha+d}}dy, \quad
II:=\int_A \frac{h(y)+h_x(y)}{|y-x|^{\alpha+d}}dy,\quad
III:=\int_E \frac{|h(y)-h_x(y)|}{|y-x|^{\alpha+d}}dy.
\end{align*}

For $I$, since $h=0$ on $B(Q, R/2)^c$, we have
\[
I\le (R/2)^{\alpha/2}
\int_{B(x, \frac{R}{8})^c} \frac{1}{|y-x|^ {\alpha+d}}dy 
+\sup_{\{z\in\bR^d:0<z_d<R/8\}}\int_{B(z, \frac{R}{8})^c\cap H^+} \frac{\delta_{H^+}(y)^{\alpha/2}}{|y-z|^{\alpha+d}}dy<\infty.
\]

For $II$, 
we first note that for any $y\in A$, 
$h(y)+h_x(y)\le c_2 |\wt{y}|^{(1+\eta)\alpha/2}$ and  $m_{d-1}(\{y:|\wt{y}|=r, -\wh \phi(\wt{y})\le y_d\le\wh \phi(\wt{y})\})\le c_3 r^{d+\eta-1}$ for $r\le R/8$, where $m_{d-1}(dy)$ is the surface measure.
Hence, for $\alpha/2<\eta$ we have
\begin{align}
{II}\le&\, c_2 \int_0^{R/8}\int_{|\wt{y}|=r} {\bf1}_A(y)|\wt{y}|^{(1+\eta)\alpha/2} |\wt{y}|^{-d-\alpha}m_{d-1}(dy)dr  \, 
\le\, c_4\int _0 ^{R/8} r^{-\alpha/2+\eta-1+\eta \alpha/2} dr<\infty.\nn
\end{align}

Next we estimate $III$.
If $0<y_d=\delta_{H^+}(y)\le \delta_D(y)$, then  $\delta_D(y)-y_d\le 4\Lambda |\wt{y}|^{1+\eta}$ and
\begin{align*}
h(y)-h_x(y)\le& (y_d+4 \Lambda|\wt{y}|^{1+\eta})^{\alpha/2}-y_d^{\alpha/2} 
\le 2 \alpha \Lambda |\wt{y}|^{1+\eta}|y_d|^{\frac{\alpha}{2} -1}.
\end{align*}
If $y_d=\delta_{H^+}(y)>\delta_D(y)$, then $\delta_D(y)\ge y_d-2\Lambda|\wt{y}|^{\eta+1}$ and
\begin{align*}
h_x(y)-h(y)\le y_d^{\alpha/2}-(y_d-2\Lambda|\wt{y}|^{\eta+1})^{\alpha/2}
\le  \alpha \Lambda|\wt{y}|^{1+\eta}(y_d-2\Lambda|\wt{y}|^{\eta+1})^{\frac{\alpha}{2}-1}.
\end{align*}
Thus, using $E\subset \{ |\wt{y}|<R/4,  \wh \phi(\wt{y})  <y_d< \wh \phi(\wt{y})+R/4\}$ and the change of variable $s=y_d-\wh \phi(r)$, we have 
\begin{align*}
III&\le c_5 \int_E \frac{|\wt{y}|^{1+\eta}(y_d-2\Lambda|\wt{y}|^{\eta+1})^{\frac{\alpha}{2} -1}}{(|\wt{y}|+|x_d-y_d|)^{d+\alpha}}dy
\le c_6 \int_0^{R/4} \int_{\wh \phi(r)}^{\wh \phi(r)+R/4}  \frac{(y_d-\wh \phi(r))^{\frac{\alpha}{2} -1}}{(r+|x_d-y_d|)^{\alpha+1-\eta}}dy_ddr\nn\\
&= c_6
\int_0^{R/4} \int_{0}^{R/4}  \frac{s^{\frac{\alpha}{2} -1}}{   (r+|x_d-(s+\wh \phi(r))|)^{\alpha+1-\eta}   }dsdr.
\end{align*}
Then, we use \cite[Lemma 4.4]{KSV8}, which is a consequence of the rearrangement inequality, and obtain 
\begin{align*}
III&\le  2 c_6
\int_0^{R/2} \left( \int_0^u t^{\frac{\alpha}{2} -1}dt \right){u^{-\alpha-1+\eta}}du 
=\frac{4 c_6}{\alpha}
\int_0^{R/2} {u^{-\frac{\alpha}{2}-1+\eta}}du <\infty.
\end{align*}
\qed

Recall that $h_{Q,r}(y)$ is defined in \eqref{e:h} for each $Q \in \partial D$ and $r \le R$. 

\begin{lemma}\label{l:sA}
For any $k >0$, let $B_k:=\left\{y \in D \cap B(Q, \frac{r}{8}): \delta_{D \cap B(Q, \frac{r}{8}) }(y) \ge 2^{-k}\right\}$.
Then, for every $|z|<2^{-k}$, 
\begin{align} \label{e:whsA}
\wh \sA_z h_{Q,r/2}(w):=\lim_{\varepsilon\to 0}\int_{|(w-z)-y|>\varepsilon} \left(h_{Q,r/2}(y)-h_{Q,r/2}(w-z)\right)\, J(w,z+ y)\, dy
\end{align}
is well defined in $B_k$.
Moreover, there exists $C_*=C_*(\eta,R,\Lambda, \rho)>0$ independent of $Q$, $k$, and $r \le R$ such that 
\[
| \wh \sA_z h_{Q,r/2}(w)| \le C_* r^{-\alpha/2} \quad \text{ for all } w \in B_k, \, |z|<2^{-k}\,.
\]
\end{lemma}

\pf
For $x \in {D}\cap B(Q, \frac{r}{8})$, let
\[
I=I(x)
:=\int_{ \R^d}{ \left|h_{Q,r/2}(y)-h_{Q,r/2}(x)\right|}\frac{|x-y|^\rho \wedge 1}{|x-y|^{d+\alpha}}  dy
\]
and
\[
II_{\varepsilon}=II_{\varepsilon} (x) 
:=\int_{ |y-x|>\varepsilon}  \left(h_{Q,r/2}(y)-h_{Q,r/2}(x)\right)\frac{dy}{|x-y|^{d+\alpha}}.
\]
For $r \le R$, let  $x^r=r^{-1} x$, $Q^r=r^{-1} Q$, and $D^r=r^{-1} D$.
The $D^r$ are $C^{1, \eta}$ open sets with the same $C^{1, \eta}$ characteristics $(1, \Lambda)$ for all $r \le R$, and 
\begin{align*}
II_{\varepsilon}  
&=r^{-\alpha/2}\int_{ |v-x^r|>\varepsilon r^{-1}} \left(\delta_{D^{r}}(v)^{\alpha/2}{\bf 1}_{D^r\cap B(Q^r, 1/2)}(v)-\delta_{D^{r}}(x^{r})^{\alpha/2} \right)\frac{ dv}{|x^r-v|^{d+\alpha} }.  
\end{align*}
Thus, by Lemma~\ref{l:frac}, $\lim_{\varepsilon \to 0}II_{\varepsilon}$ exists and satisfies $|\lim_{\varepsilon \to 0}II_{\varepsilon}| \le c_1 r^{-\alpha/2}$.

Similarly, we obtain  
\begin{align*}
I
&=r^{-\alpha/2}\int_{ \R^d}\left|\delta_{D^{r}}(v)^{\alpha/2}{\bf 1}_{D^r\cap B(Q^r, 1/2)}(v)-\delta_{D^{r}}(x^{r})^{\alpha/2}\right| \frac{r^\rho|x^r-v|^\rho \wedge 1}{|x^r-v|^{d+\alpha}}  dv.
\end{align*}
Since $|\delta_{D^{r}}(v)^{\alpha/2}-\delta_{D^{r}}(x^r)^{\alpha/2}|\le |\delta_{D^{r}}(v)-\delta_{D^{r}}(x^r)|^{\alpha/2}\le |v-x^r|^{\alpha/2}$, for $r\le 1$ we have 
\begin{align*}
&\int_{ \R^d}\left|\delta_{D^{r}}(v)^{\alpha/2}{\bf 1}_{D^r\cap B(Q^r, 1/2)}(v)-\delta_{D^{r}}(x^r)^{\alpha/2}\right| \frac{r^\rho|x^r-v|^\rho \wedge 1}{|x^r-v|^{d+\alpha}}  dv\\
&\le \int_{ D^r\cap B(Q^r, 1/2) }\frac{|x^r-v|^{\rho+\alpha/2}}{|x^r-v|^{d+\alpha}}  dv+c_2\int_{(D^r\cap B(Q^r, 1/2))^c } \,\frac{1}{|x^r-v|^{d+\alpha}}  dv\\
& \le c_3\, \int_{ \R^d}\frac{|u|^{\rho+\alpha/2} \wedge 1}{|u|^{d+\alpha}}  du <\infty.
\end{align*}
In the last inequality above we used the assumption that $\rho>\alpha/2$. 

From \eqref{eqn:exp}--\eqref{e:J2} we observe that
\begin{align*}
&\int_{\{y:|(w-z)-y|>\varepsilon\}} \left(h_{Q,r/2}(y)-h_{Q,r/2}(w-z)\right)\, J(w,z+ y)\, dy\\
=&\int_{\{y:|(w-z)-y|>\varepsilon\}} \left(h_{Q,r/2}(y)-h_{Q,r/2}(w-z)\right)\, 
\frac{\kappa(w,z+y) }{ |w-z-y|^{d+\alpha} \psi_1 (|w-z-y|)}dy
 \\
= & \int_{\{y:|(w-z)-y|>\varepsilon\}} \left(\kappa(w, z+y)
-\kappa(w,w)\right)\frac{(h_{Q,r/2}(y)-h_{Q,r/2}(w-z))}{\psi_1(|w-z-y|) |w-z-y|^{d+\alpha}}dy \\
&\qquad+ \, \kappa(w,w)
\int_{\{y:|(w-z)-y|>\varepsilon\}}
\frac{(h_{Q,r/2}(y)-h_{Q,r/2}(w-z))}{\psi_1(w-z-y)|w-z-y|^{d+\alpha}}dy
\end{align*}
and
\[
\left|\int_{\{y:|(w-z)-y|>\varepsilon\}} \left(h_{Q,r/2}(y)-h_{Q,r/2}(w-z)\right)\, J(w,z+ y)\, dy\right|
\le c_4 I(w-z)+c_4  II_{\varepsilon} (w-z).
\]
Therefore, $\wh\sA_z h_{Q,r/2} (w)$ exists 
on $B_k$
and we have $ |\wh \sA_z h_{Q,r/2} (w)|\le  c_5 r^{-\alpha/2}$ for every $w \in B_k$ and  $|z|<2^{-k}$.
\qed

Using Lemma~\ref{l:sA}, we prove the following theorem which plays a critical role in estimating the exit distribution. 
In the next theorem for $x\in D$, we use $z_x$ to denote a point on $ \partial D$ such that  $|z_x-x|=\delta_D(x)$,
and we use the  coordinate system $CS_{z_x}$ with a $C^{1,\eta}$ function $\phi$ such that $\phi(0)=0$, $ \nabla\phi (0)=(0, \dots, 0)$, $\| \nabla\phi \|_\infty \leq \Lambda$, $| \nabla \phi (\wt y)-\nabla \phi (\wt w)| \leq \Lambda |\wt y-\wt w|^{\eta}$, and $B(0, R)\cap {D}= \{ y=(\wt y, \, y_d) \in B(z_x, R) \mbox{ in } CS_{z_x} : \phi( \wt y )< y_d  \}$. 
For the next theorem and its proof, we always use this coordinate system $CS_{z_x}$. 

\begin{thm}\label{L:2}
There are constants 
$b_1=b_1(\eta,R,\Lambda, \rho) \in (0, 1/10)$ and $c_1=c_1(\eta,R,\Lambda)>1$ 
such that for any $r \le b_1(R \wedge 1)/2$ and $x \in D$ with $\delta_D(x)<r$ we have
\begin{equation}\label{e:L:3}
\E_x\left[ \tau_{D \cap B(z_x, r)}\right] \le c_1 \, r^{\alpha/2}\delta_{D}(x)^{\alpha/2}\qquad \hbox{where } z_x \in \partial D  \text{ with }
\delta_D(x)=|x-z_x|,
\end{equation}
and for any $r\le (R \wedge 1)/4$, $\lambda \ge 4$ and $x \in D$ with $\delta_D(x)<\lambda^{-1} r/2$ we have 
\begin{align}\label{e:L:3-1}
\P_x\left(Y_{\tau_{D \cap B(z_x, \lambda^{-1}r)}} \in 
\{2 \Lambda |\wt y|<y_d, \lambda^{-1}r <|y|<r  \mbox{ in } CS_{z_x}\}\right) \ge&  c_1^{-1}\frac{\delta_{D}(x)^{\alpha/2}}{r^{\alpha/2}}\, 
\end{align}
where $z_x \in \partial D$  and $\delta_D(x)=|x-z_x|$.  
\end{thm}
\pf
Without loss of generality, we assume that $z_x=0$ and let $A(a,b):=B(0, b)\backslash B(0,a)$ with $0<a<b$.
Let $r \le (R \wedge 1)/2$
 and $h(y)=h_{0,r/2}(y)$ (see \eqref{e:h}).
Let $f$ be a nonnegative smooth radial function such that $f(y)=0$ for $|y|>1$ and $\int_{\R^d} f (y) dy=1$.
For $k\geq 1$, we define $f_k(y):=2^{kd} f (2^k y)$ and $h^{(k)}(z):= ( f_k*h)(z)\in C_c^\infty(\R^d)$, and we let $B_k:=\{y \in D \cap B(0, {r}/{8}): \delta_{D \cap B(0, {r}/{8}) }(y) \ge 2^{-k}\}$.

By Lemma~\ref{l:sA},  $\wh \sA_z h(w)$ exists for $w\in B_k$ and $z\in B(0,2^{-k})$, with $-C_* r^{-\alpha/2}\le \wh \sA_z h(w) \le C_* r^{-\alpha/2}$, where $\wh \sA_z h(w)$ is defined in \eqref{e:whsA} and  $C_{*}$ is the constant in  Lemma~\ref{l:sA}.
Then, by letting $\varepsilon \to 0$ and using the dominated convergence theorem, it follows that $\sA h^{(k)}$ is well defined everywhere and for large $k$ and $|z|<2^{-k}$ we have
\begin{equation} \label{e:*333}
| \sA h^{(k)}(w)|=|\int_{|z|<2^{-k}} f_k(z)\wh \sA_z h(w)\, dz|\le
C_* r^{-\alpha/2} \int_{|z|<2^{-k}} f_k(z) \, dz \le 
C_* r^{-\alpha/2} \quad \text{ on } B_k\, .
\end{equation}
Applying \eqref{e:334} to $U^k_\lambda :=D \cap B(0, \lambda^{-1}r)\cap B_k$ with $\lambda \ge 8$ and $h^{(k)}$ and using \eqref{e:*333} we have
\begin{equation*}\label{e:sh}
\E_{x}\left[h^{(k)}\big(Y_{\tau_{U^k_\lambda}}\big)\right]-C_{*} r^{-\alpha/2} \E_x\left[\tau_{U^k_\lambda}\right] \le h^{(k)}(x) \le \E_{x}\left[h^{(k)}\big(Y_{\tau_{U^k_\lambda}}\big)\right]+C_{*} r^{-\alpha/2} 
\E_x\left[\tau_{U^k_\lambda}\right], \quad x \in U^k_\lambda.
\end{equation*}
Since $h^{(k)}$ is in $C^\infty_c(\R^d)$ and by letting $k\to \infty$, for all $\lambda \ge 8$ and $x \in D \cap B(0, \lambda^{-1}r)$ we obtain 
\begin{align}
    &\delta_{D}(x)^{\alpha/2}  \ge \E_{x}\left[h\left(Y_{ \tau_{D \cap B(0, \lambda^{-1}r)}} \right)  \right] - C_{*}  r^{-\alpha/2} \E_x\left[\tau_{D \cap B(0, \lambda^{-1}r)}\right] \label{e:ss}
\end{align}
and
\begin{align}
    &\delta_{D}(x)^{\alpha/2} \le   \E_{x}\left[h\left(Y_{ \tau_{D \cap B(0, \lambda^{-1}r)}} \right)  \right]+ C_{*}  r^{-\alpha/2}\E_x\left[\tau_{D \cap B(0, \lambda^{-1}r)}\right]    \label{e:ss1}.
\end{align}

For any $z  \in D \cap B(0, \lambda^{-1}r)$ and $y \in D\cap (B(0, 2^{-1} r)\backslash B(0,  \lambda^{-1} r))$,  
since $2|y| \le r \le 1/2$,  we have $j(|y-z| ) \ge j(|y|+|z|) \ge   j(2|y|) \ge  c_1 j(|y|)$.
Thus, by \eqref{e:levy} we obtain 
\begin{align} \E_x\left[h\left(Y_{\tau_{D\cap B(0, \lambda^{-1}r)}}\right)\right]
&= \E_x \int_{ D\cap A(\lambda^{-1}r, 2^{-1} r)} \int_0^{\tau_{D \cap B(0, \lambda^{-1}r)}}  j(|Y_t-y|)dt \delta_{D}(y)^{\alpha/2} dy \nonumber \\
&\ge  c_1  \, \E_x\left[\tau_{D \cap B(0, \lambda^{-1}r)}\right] \int_{  D\cap A(\lambda^{-1}r, 2^{-1} r)}    j(|y|)\delta_{D}(y)^{\alpha/2} dy. \label{e:new23}
\end{align}
Similarly,  with $V:=\{2 \Lambda |\wt y|<y_d\}$ we also have
\begin{align}
\P_x\left(Y_{\tau_{D\cap B(0, \lambda^{-1}r)}}\in V \cap A(\lambda^{-1}r, 2^{-1} r)
\right)
 \ge  c_2  \,\E_x\left[\tau_{D \cap B(0, \lambda^{-1}r)}\right] \int_{ V\cap A(\lambda^{-1}r, 2^{-1} r)}  j(|y|) dy. \label{e:new23-1}
\end{align}

Clearly, 
\begin{align}
 \int_{ V\cap A(\lambda^{-1}r, 2^{-1} r)} |y|^{-d-\alpha} dy
\ge c_3r^{-\alpha} \left(\lambda^{\alpha}-1\right)
\, . \label{e:new25-1}
\end{align}
Since for every $y\in B(0, R) \cap D$ with $ 2 \Lambda |\wt y|<y_d$ we have 
\[
\delta_D(y)\ge (1+\Lambda)^{-1}\left(y_d-\phi(\wt y)\right) 
\ge  (2\Lambda)^{-1} (y_d -\Lambda|\wt y|)> (4\Lambda)^{-1}|y_d|\ge  ((2\Lambda)^{-2}+1)^{-1/2}(4\Lambda)^{-1}|y|, 
\]
by changing to polar coordinates with $|y|=s$ we obtain
\begin{align}
&\int_{  D\cap A(\lambda^{-1}r, 2^{-1} r)}    j(|y|)(\delta_{D}(y))^{\alpha/2} dy \ge\int_{ \{(\wt y, y_d)\in D: 2 \Lambda|\wt y|< y_d,   \lambda^{-1} r<|y| <2^{-1} r\}} j(|y|)\delta_{D}(y)^{\alpha/2}  dy\nn\\
\ge& c_4\int_{\{(\wt y, y_d): 2\Lambda |\wt y|< y_d,   \lambda^{-1} r<|y| <2^{-1} r\}} |y|^{-d-\alpha} |y|^{\alpha/2} dy
 \ge c_ 5r^{-\alpha/2} \left(\lambda^{\alpha/2}-1\right)\, \label{e:new25}.
\end{align}
Then, combining \eqref{e:new23} and \eqref{e:new25} yields
\begin{align}
  \E_{x}\left[h\left(Y_{ \tau_{D \cap B(0, \lambda^{-1}r)}} \right) \right]\ge
    c_{6} r^{-\alpha/2}\left(\lambda^{\alpha/2}-1\right)
\E_x\left[\tau_{D \cap B(0, \lambda^{-1}r)}\right]\, \label{e:new26}
\end{align}
and \eqref{e:new23-1} and \eqref{e:new25-1} yield
\begin{align}
\P_x&\left(Y_{\tau_{D\cap B(0, \lambda^{-1}R)}} \in   V\cap A(\lambda^{-1}r, 2^{-1} r)\right) 
\ge c_6 r^{-\alpha} \left(\lambda^{\alpha}-1\right)
\E_x\left[\tau_{D \cap B(0, \lambda^{-1}r)}\right].  \label{e:new26-1}
\end{align}
Hence, by \eqref{e:ss} and \eqref{e:new26}, we find that for every $\lambda \ge 
\lambda_0:=(2 + 2C_*/c_6 )^{2/\alpha} \vee (10)$ 
and $\delta_D(x)<\lambda^{-1}r$ we have
\begin{align}
\delta_{D}(x)^{\alpha/2}\ge 
\left(c_6 \lambda^{\alpha/2} -(c_6+C_*) \right)
r^{-\alpha/2} \E_x[\tau_{D \cap B(0, \lambda^{-1}r)}] \ge 
\frac{c_6}{2} (\lambda^{-1} r)^{-\alpha/2}\E_x[\tau_{D \cap B(0, \lambda^{-1}r)}] .
\label{e:ss6}
\end{align}
Thus, we have proved \eqref{e:L:3} with $b_1=\lambda_0^{-1}$ and $r=(R \wedge 1)/2$.

Conversely, using \eqref{e:ss1} and then using Lemma~\ref{L:2.3} and \eqref{e:new26-1}, we find that for every $\lambda \ge 8$ and $\delta_D(x)<\lambda^{-1}r/2$ we obtain
\begin{align*}
\delta_D(x)^{\alpha/2}
\le& (r/2)^{\alpha/2}
\P_x\left(Y_{\tau_{D \cap B(0, \lambda^{-1}r)}} \in  D\right) +C_{*} r^{-\alpha/2}\E_x\left[\tau_{D \cap B(0, \lambda^{-1}r)}\right]\nn\\
\le&  r^{-\alpha/2} \left(c_7\lambda^{\alpha} +C_{*}  \right)\E_x\left[\tau_{D \cap B(0, \lambda^{-1}r)}\right]\nn\\
\le &r^{\alpha/2}\frac{c_7\lambda^{\alpha} +C_*}
{(\lambda^\alpha-1)c_{6}} 
\P_x\left(Y_{\tau_{D \cap B(0, \lambda^{-1}r)}} \in  V\cap A(\lambda^{-1}r, 2^{-1} r)\right)\nn \\
\le &r^{\alpha/2}\frac{(c_7+C_*)\lambda^{\alpha}}
{(\lambda^\alpha-(\lambda/2)^\alpha)c_{6}}
\P_x\left(Y_{\tau_{D \cap B(0, \lambda^{-1}r)}} \in  V\cap A(\lambda^{-1}r, 2^{-1} r)\right)\nn \\
= &r^{\alpha/2}
\frac{c_7+C_*}{(1-2^{-\alpha})c_{6}}
\P_x\left(Y_{\tau_{D \cap B(0, \lambda^{-1}r)}} \in  V\cap A(\lambda^{-1}r, 2^{-1} r)\right).
\end{align*}
Thus, we have proved \eqref{e:L:3-1}.
\qed

\section{Preliminary lower bound estimates }\label{s:plbd}

In this section, we discuss a preliminary lower bound for $p_D(t, x, y)$.

Using \cite[Theorem 1.4 and Lemma 2.5]{CKK2}, the proof of the next lemma is the same as that of \cite[Lemma 3.1]{CKS2}. Thus, we omit the proof. 

\begin{lemma}\label{L:3.3}
Let $T$, $a$, and $b$ be positive constants.
For any $\b\in[0,\infty]$, there exists a constant $c=c(a, b, \b, T)>0$ such that for all $\lambda \in (0, T]$ we have
\[
\inf_{y\in\bR^d \atop |y -z| \le b \lambda^{1/\alpha}} \P_y \left(\tau_{B(z, 2b \lambda^{1/\alpha} )} > a\lambda \right)\, \ge\, c_.
\]
\end{lemma}

Next, we give some preliminary lower bound estimates for $p_D(t, x, y)$ on $\delta_D(x)\wedge \delta_D(y)\wedge T\ge t^{1/\alpha}$, which are used to derive the sharp two-sided estimates for $p_D(t, x, y)$. 
We first consider $D$ an arbitrary nonempty open set, and we use the convention that $\delta_{D}(\cdot) \equiv \infty$ when $D =\bR^d$.
This convention allows us to derive the lower bound of $p(t,x,y)$ simultaneously. 

Using \cite[Theorem 1.4]{CKK2} and Lemma~\ref{L:3.3}, the proof of the next lemma is the same as that of \cite[Proposition 3.2]{CKS2}.
Thus, we omit the proof. 

\begin{prop}\label{step1}
Let $D$ be an arbitrary open set and let $a$ and $T$ be positive constants.
Suppose that $(t, x, y)\in (0, T]\times D\times D$, with $\delta_D(x) \ge a t^{1/\alpha} \geq 2|x-y|$.
Then, for any $\b \in [0, \infty]$, there exists a positive constant $c=c(a,  \b, T)$ such that $p_D(t,x,y) \,\ge\,c\, t^{-d/\alpha}$.
\end{prop}

\begin{prop}\label{step3}
Let $D$ be an arbitrary open set and let $a$ and $T$ be positive constants. 
Suppose that $(t, x, y)\in (0, T]\times D\times D$, with $\delta_D(x)\wedge \delta_D (y) \ge a t^{1/\alpha}$ and $a t^{1/\alpha} \leq 2|x-y|$.
Then, for any $\b \in [0, \infty]$, there exists a constant $c=c(a, \b, T)>0$ such that $p_D(t, x, y)\ge c t j(|x-y|)$.
\end{prop}

\pf 
By Lemma~\ref{L:3.3}, starting at $z\in B(y, \, 4^{-1} a t^{1/\alpha})$, with probability at least $c_1=c_1(a, \b, T)>0$ the process $Y$ does not move more than $6^{-1} a t^{1/\alpha} $ by time $t$.
Thus, using the strong Markov property and the L\'evy system in \eqref{e:levy}, we obtain
\begin{align}
&\P_x \left( Y^D_t \in B \big( y, \,  2^{-1} a t^{1/\alpha} \big) \right)\nn\\  
&\ge  c_1\P_x(Y_{t\wedge \tau_{B(x, 6^{-1} a t^{1/\alpha})}}^D\in B(y,\, 4^{-1}a t^{1/\alpha})\hbox{ and }t \wedge \tau_{B(x, 6^{-1}a  t^{1/\alpha})} \hbox{ is a jumping time })\nonumber \\
&= c_1 \E_x \left[\int_0^{t\wedge \tau_{B(x, 6^{-1}a t^{1/\alpha})}} \int_{B(y, \, 4^{-1}a t^{1/\alpha})}
J(Y_s, u) duds \right]. \label{e:nv1}
\end{align}
Lemma~\ref{L:3.3} also implies that 
\begin{equation}\label{eq:lowtau}
\E_x \left[ t \wedge
\tau_{B(x, 6^{-1}a t^{1/\alpha})} \right] \,\ge\,t \, \P_x
\left(\tau_{B(x, 6^{-1}a t^{1/\alpha})} \ge  t \right)
 \,\ge\, c_2\,t
\qquad \hbox{ for all } t\in (0, T].
\end{equation}

We fix the point $w$ on the line connecting $|x-y|$ (i.e., $|x-y|=|x-w|+|w-y|$) such that $|w-y|=7 \cdot 2^{-5} a t^{1/\alpha}$, which is possible because  $\delta_D (y) \ge a t^{1/\alpha}$.
Then, $B(w, 2^{-5} a t^{1/\alpha}) \subset B(y, \, 4^{-1}a t^{1/\alpha})$.
Moreover, for every $(z,u) \in  B(x, 6^{-1}a t^{1/\alpha}) \times B(w, 2^{-5} a t^{1/\alpha}) $ we have 
\[
|z-u| < 6^{-1}a t^{1/\alpha}+ 2^{-5} a t^{1/\alpha}+|x-w| =|x-y|+(6^{-1}+2^{-5}-7 \cdot 2^{-5})a t^{1/\alpha} <|x-y|.
\]
Thus, $B(w, 2^{-5} a t^{1/\alpha}) \subset B(y, \, 4^{-1}a t^{1/\alpha}) \cap \{u:|u-z|<|x-y|\}$.
Combining this result with \eqref{e:conkappa1} and \eqref{eq:lowtau}, we obtain 
\begin{align}
 &\E_x \left[\int_0^{t\wedge \tau_{B(x, 6^{-1}a t^{1/\alpha})}} \int_{B(y, \, 4^{-1}a t^{1/\alpha})}
J(Y_s, u) duds \right] \nn\\
&\ge \E_x \left[\int_0^{t\wedge \tau_{B(x, 6^{-1}a t^{1/\alpha})}} \int_{B(w, 2^{-5} a t^{1/\alpha})}
J(Y_s, u) {\bf 1}_{\{ |Y_s -u| <|x-y| \}} duds \right] \nn\\
&\ge c_3 \E_x \left[ t \wedge
\tau_{B(x, 6^{-1}a t^{1/\alpha})} \right]   |B(w, 2^{-5} a t^{1/\alpha})| j(|x-y|) >c_4    t^{1+d/\alpha} j(|x-y|) .\label{e:nv2}
\end{align} 
Then, using the semigroup property along with \eqref{e:nv2} and Proposition~\ref{step1}, the proposition follows from the proof of \cite[Proposition 3.4]{CKS2}.
\qed

Combining Propositions~\ref{step1} and~\ref{step3} with the definition of $j$, we obtain a lower bound for $p_D(t,x,y)$ that yields the preliminary lower bound for $p_D(t,x,y)$ and $p(t, x, y)$ for the case $\b\in[0, 1]$ and the case $\b\in(1, \infty]$ with $|x-y|<1$. 

\begin{prop}\label{step4}
Let $D$ be an arbitrary open set and let $a$ and $T$ be positive constants.
Suppose that $(t, x, y)\in (0, T]\times D\times D$, with $\delta_D(x)\wedge \delta_D(y) \ge a t^{1/\alpha}$.
Then, for any $\b\in[0,\infty]$, there exists a positive constant $c=c(a, \beta, T)$ such that
\begin{align*}
p_D(t, x, y)\,\ge \,
c \left(  t^{-d/\alpha}\wedge t j(|x-y|) \right).
\end{align*}
\end{prop}

We next consider cases $\b\in(1, \infty]$ with  $|x-y|\ge 1$.
We will closely follow the proofs of \cite[Theorem 3.6]{CKK} and \cite[Theorem 5.5]{CKK3}.

For the remainder of this section, we assume that $D$ is an open set with the following property: there exist $\lambda_1 \in [1, \infty)$ and $\lambda_2 \in (0, 1]$ such that for every $r \le 1$ and  $x,y$ in the same component of $D$ with $\delta_D(x)\wedge \delta_D(y)\ge r$ there exists in $D$ a length parameterized rectifiable curve $l$ connecting $x$ to $y$ with the length $|l|$ of $l$ is less than or equal to $\lambda_1|x-y|$ and $\delta_D(l(u))\geq\lambda_2 r$ for $u\in[0,|l|].$ 

Under this assumption, we prove the preliminary lower bound of $p_D(t,x,y)$ on $|x-y|\ge 1$ separately for the case $\b=\infty$ and the case $\b\in(1, \infty)$. 

\begin{prop}\label{step5}
Suppose that $T>0$, $a \in (0, 4^{-1}T^{-1/\alpha}]$, and $\b=\infty$. 
Then, there exist constants $c_i=c_i(a, T, \lambda_1, \lambda_2)>0$, $i=1,2$, such that for any $x, y$ in the same component of $D$ with $\delta_D(x)\wedge\delta_D(y) \ge a t^{1/\alpha}$, $|x-y|\ge 1$, and $t \le T$ we have  
\[
p_D(t, x, y)\,\ge \,c_{1} \left(\frac{t}{T|x-y|}\right)^{c_{2}|x-y|}.
\]
\end{prop}
\pf 
We fix $T>0$ and $a \in (0, 4^{-1}T^{-1/\alpha}]$, and we let $R_1:=|x-y| \ge 1$.
By our assumption for $D$, there is a length parameterized curve $l\subset D$ connecting $x$ and $y$ such that the total length $|l|$ of $l$ is less than or equal to $\lambda_1 R_1$ and $\delta_D(l(u))\ge \lambda_2 at^{1/\alpha}$ for every $u\in [0, |l|]$.
We define $k$ as the integer satisfying $(4 \le ) 4 \lambda_1 R_1\leq k <4 \lambda_1 R_1+1<5\lambda_1 R_1$ and $r_t:= 2^{-1}\lambda_2 a t^{1/\alpha} \le 8^{-1}$.
Let $x_i:=l(i|l|/k)$ and  $B_i:=B(x_i, r_t)$, with $i=0,1,2,\ldots,k$.
Then, $\delta_D(x_i)> 2r_t$ and $B_i=B(x_i, r_t) \subset B(x_i, 2 r_t)\subset D$, with $i=0,1,2,\ldots,k$. 

Since $4 \lambda_1 R_1\le k$, for each $y_i\in B_i$ we have  
\begin{eqnarray}\label{e:stst1}
|y_i-y_{i+1}|\leq |y_i-x_i|+|x_i-x_{i+1}|+|x_{i+1}-y_{i+1}| \leq 
\frac{1}{8}+\frac{|l|}{k}+ \frac{1}{8}<\frac{\lambda_1 R_1}{4\lambda_1 R_1}+\frac{1}{4}\le\frac{1}{2} .
\end{eqnarray}
Moreover, $\delta_D(y_i)\ge\delta_D(x_i)-|y_i-x_i |>r_t>r_{t/k}$.

Thus, by Proposition~\ref{step4} and~\eqref{e:stst1}, there are constants $c_i=c_i(a, T,  \lambda_2)>0$, $i=1,2$, such that for $(y_i, y_{i+1})\in B_i\times B_{i+1}$ we have
\begin{eqnarray}\label{e:stst2}
p_D(t/k, y_i, y_{i+1})\ge c_1\left((t/k)^{-d/\alpha} \wedge \frac{t/k}{|y_i-y_{i+1}|^{d+\alpha}}\right)\ge  c_2 \,t/(Tk).
\end{eqnarray}
Observe that $4\lambda_1R_1\le k < 2(k-1) < 8\lambda_1R_1$ and $r_t\ge T^{1/\alpha} r_{t/(Tk)}$. 
Thus, from \eqref{e:stst2} we obtain
\begin{align*}
&p_D(t,x,y)\ge \int_{B_1}\ldots\int_{B_{k-1}}p_D(t/k,x,y_1)\ldots p_D(t/k, y_{k-1},y) dy_{k-1}\ldots dy_1\\
&\ge (c_2t(Tk)^{-1})^k\Pi^{k-1}_{i=1} |B_i| \ge  c_3 (c_4t(Tk)^{-1})^{c_5k}  \ge  c_6 (c_7t(TR_1)^{-1})^{c_8R_1}  \ge c_9 (t(TR_1)^{-1})^{c_{10}R_1}.
\end{align*}
\qed

\begin{prop}\label{step6}
Suppose that $T>0$, $a \in (0, 4^{-1}T^{-1/\alpha}]$, and $\beta\in(1, \infty)$. 
Then, there exist constants  $c_i=c_i(a, \beta,T,  \lambda_1, \lambda_2)>0$, $i=1,2$ such that for any $x, y$ in the same component of $D$ with $\delta_D(x)\wedge\delta_D(y) \ge a t^{1/\alpha}$, $|x-y|\ge 1$, and $t \le T$ we have  
\begin{eqnarray*}
p_D(t,x,y)\ge c_1 t \exp\left\{ -c_2 \left(|x-y|\left(\log\frac{T|x-y|}{t}\right)^{\frac{\beta -1}{\beta}}\wedge (|x-y|)^{\beta}\right)\right\}.
\end{eqnarray*}
\end{prop} 

\pf
We fix $T>0$ and $a \in (0, 4^{-1}T^{-1/\alpha}]$, and we let $R_1:=|x-y|$. 
If either $1 \le R_1 \le 2$ or $R_1(\log(TR_1/t))^{(\beta-1)/\beta}\ge (R_1)^{\beta}$, the proposition holds by virtue of Proposition~\ref{step4}. 
Thus, for the remainder of this proof we assume that $R_1> 2$ and $R_1(\log(TR_1/t))^{(\beta-1)/\beta}< (R_1)^{\beta}$, which is equivalent to $R_1 \exp\{-\left(R_1\right)^{\beta}\}< t/T$.

Let $k\ge 2$ be a positive integer such that
\begin{align}\label{e:k}
1<R_1 \left(\log\frac{TR_1}{t}\right)^{-1/\beta}\le k<R_1 \left(\log\frac{TR_1}{t}\right)^{-1/\beta}+1< 2R_1 \left(\log\frac{TR_1}{t}\right)^{-1/\beta}.
\end{align} 
By our assumption for $D$, there is a length parameterized curve $l\subset D$ connecting $x$ and $y$ such that the total length $|l|$ of $l$ is less than or equal to $\lambda_1 R_1$ and $\delta_D(l(u))\ge \lambda_2 at^{1/\alpha}$ for every $u\in [0, |l|]$.
We define 
$r_t:=(2^{-1}{\lambda_2}at^{1/\alpha})\wedge ((6\lambda_1)^{-1}(\log ({TR_1}/{t}))^{1/\b})$.
Then, by \eqref{e:k} and the assumption $((\log(TR_1/t))^{1/\beta}) \vee 2 < R_1$ we have
\begin{align}
\left(\frac{\lambda_2}{2} a{T}^{1/\alpha}
\left(\frac{t}{TR_1}\right)^{1/\alpha}\right)\wedge \left(\frac{(2 \log 2)^{1/\b}}{6\lambda_1}\left( \frac{t}{TR_1}\right)^{1/\b} \right) \le r_t\le\frac{1}{6\lambda_1}\left(\log\frac{TR_1}{t} \right)^{1/\b}< \frac{R_1}{3\lambda_1 k}.\label{e:r_t.2}
\end{align}

We define $x_i:=l(i |l|/k)$ and $B_i:=B(x_i, r_t)$, with $i=0, \ldots, k$.
Then, $\delta_D(y_i)\ge 2^{-1} \lambda_2 a t^{1/\alpha}>2^{-1} \lambda_2 a (t/k)^{1/\alpha} $ for every  $y_i \in B_i$. 
Note that from \eqref{e:r_t.2} we obtain 
\begin{align}\label{e:y_i}
|y_i-y_{i+1}|\le|x_i-x_{i+1}|+2r_t \le\left(\lambda_1 +\frac{2}{3 \lambda_1} \right)\frac{R_1}{k}. 
\end{align}

Thus, using Proposition~\ref{step4} along with \eqref{e:k} and  \eqref{e:y_i} we obtain
\begin{align}\label{p(t/k)}
&p_D(t/k, y_i, y_{i+1}) 
\ge c_{1}\left((t/k)^{-d/\alpha}\wedge\frac{t}{k} 
j(|y_i-y_{i+1}|)
\right)\ge \, c_{2} \left(1 \wedge (\frac{t}{k}\left(R_1/k\right)^{-d-\alpha}e^{-c_{3}(R_1/k)^{\b}})\right)\nn\\
& \ge \,c_{4} \frac{t}{TR_1}
\left(\frac{k}{R_1}\right)^{d+\alpha-1}
e^{-c_{3}(R_1/k)^{\b}}\, \ge \, c_{4} \frac{t}{TR_1}\left(\log \frac{TR_1}{t}\right)^{-\frac{d+\alpha-1}{\b}}\left(\frac{t}{TR_1}\right)^{c_{3}}\ge \, c_{4}\left(\frac{t}{TR_1}\right)^{c_{5}}.
\end{align}
Since the lower bound of $r_t$ in \eqref{e:r_t.2} yields $r_t \ge c_{6}(t/(TR_1))^{(\alpha \wedge \beta)^{-1}}$, by using \eqref{p(t/k)} and the semigroup property we conclude that 
\begin{eqnarray*}
p_D(t, x, y)&\ge& \int_{B_1}\cdots\int_{B_{k-1}} p_D(t/k, x, y_1)\cdots p_D(t/k, y_{k-1}, y ) dy_1\cdots dy_{k-1}\\
&\ge& c_{7}\exp\{-c_{8}k\log({TR_1/t})\}\\
&\ge& c_{7}\exp\left\{-c_{8}\left(R_1 \log\left(\frac{TR_1}{t}\right)^{-1/\beta}+1\right)\log\frac{TR_1}{t}\right\}\\
&\ge&c_{7} \exp\left\{-c_{9}\left(R_1 \log\left(\frac{TR_1}{t}\right)^{1-1/\beta}\right)\right\}.
\end{eqnarray*}
\qed

\medskip

{\bf Proof of the lower bound in \eqref{eq:M-THM1}}.
The proof for the  two cases $\b\in[0, 1]$ and $\b\in(1, \infty]$ with $|x-y|<1$  follow from Proposition~\ref{step4} with $D=\R^d$. 
The proof for the remaining cases  follows from Propositions~\ref{step5} and~\ref{step6} with $D=\R^d$. 
\qed

\section{Upper bound estimates}

In this section, we derive the upper bound estimate for $p_D(t, x, y)$ as stated in Theorem~\ref{t:main}.
We first introduce a lemma that appears in~\cite{CKS6}.
The proof of the next lemma is identical to that of \cite[Lemma 3.1]{CKS6}, so we omit the proof. 

\begin{lemma}\label{L:4.1} 
Suppose that $U_1, U_3, E$ are open subsets of $\mathbb{R}^d$, with $U_1, U_3\subset E$ and $dist (U_1, U_3)>0$.
Let $U_2:=E\backslash(U_1\cup U_3)$.
If $x\in U_1$ and $y\in U_3$, then for every $t>0$ we have
\begin{align}
p_E(t, x, y) &\le \mathbb{P}_x \left(Y_{\tau_{U_1}}\in U_2\right)\cdot\sup _{s<t, z\in U_2} p_E(s, z, y)\nn\\
&+ \int_0^t \mathbb{P}_x\left(\tau_{U_1}>s\right) \mathbb{P}_y \left(\tau_E>t-s \right)ds \cdot\sup_{u\in U_1, z\in U_3} J(u, z) \label{eq:ub1}\\
&\le \mathbb{P}_x\left(Y_{\tau_{U_1}}\in U_2\right) \cdot\sup_{s<t, z\in U_2} p(s, z, y)+\left(t \wedge \mathbb{E}_x[\tau_{U_1}]\right)\cdot\sup_{u\in U_1, z\in U_3} J(u, z). \label{eq:ub2}
\end{align}
\end{lemma}

For the remainder of this section we assume that $\eta\in (\alpha/2, 1]$, $T>0$, and  $D$ is a $C^{1,\eta}$ open set with characteristics $(R, \Lambda)$.
Without loss of generality, we assume that $\Lambda>1$ and $R <10^{-1}$.
Recall that $b_1$ is the constant in Theorem~\ref{L:2}.
We let 
\begin{equation*}
a_T=a_{T,R}:=2^{-1} b_1 {R}{T^{-1/\alpha}} <(200)^{-1}T^{-1/\alpha},
\end{equation*}
and for $x\in D$ we use $z_x$ to denote a point on $\partial D$ such that  $|z_x-x|=\delta_D(x)$.

We first obtain the upper bound for the survival probability.
Recall that $\Psi$ is defined in~\eqref{e:dax}.
\begin{lemma}\label{l:sdfnW}
There exists a positive constant $c=c(\beta, R, \Lambda, \eta, \rho, T)$ such that for any $(t, x)\in (0, T]\times D$ we have
$\P_x(\tau_D >t)  \le c \Psi(t,x).$
\end{lemma}
\pf
We need to prove the lemma only for $\delta_D(x)  \le a_Tt^{1/\alpha}/8$. 
Let 
$U:=D\cap B(z_x, a_Tt^{1/\alpha})$.
Since $\P_x(\tau_D >t) \le \P_x(\tau_U >t)+\P_x(X_{\tau_U}  \in D)$, by Chebyshev's inequality, Lemma~\ref{L:2.3}, and \eqref{e:L:3} we have  
$
\P_x(\tau_D >t) 
\le   t^{-1}\E_x[\tau_U]  + c_1(a_Tt^{1/\alpha})^{-\alpha }\E_x[\tau_U]  \le  c_2   \delta_D(x)^{\alpha/2}/\sqrt{t} \le c_3 \Psi(t,x). 
$
\qed

Next, we use \eqref{eq:ub2} to obtain the intermediate upper bound in which one boundary decay appears. 

\begin{prop}\label{Pr:4.4}
For any  $a\le a_T$ and  $\b\in[0,\infty]$, there exists a positive constant $c=c(\b, R, \Lambda, T, \eta, \rho, a)$ such that for every $(t, x, y)\in (0, T]\times D\times D$ with $|x-y|\ge  12at^{1/\alpha}{\bf 1}_{\b\in[0, 1]} +2 \cdot {\bf 1}_{\b\in(1, \infty)}+2(1+2at^{1/\alpha} )\cdot {\bf 1}_{\b= \infty}$ we have 
\[
p_D(t, x, y)\le c \Psi(t,x)\cdot 
\begin{cases}
h_{ C_1\wedge\gamma_1,\gamma_1, T}(t, |x-y|/3)&\mbox{ if } \b\in[0,\infty),\\
h_{ C_1, \gamma_1, T}(t, |x-y|/2)&\mbox{ if } \b=\infty,
\end{cases}
\]
where $C_1$ is the constant in Theorem~\ref{T:1.1} and $\gamma_1$ is the constant in~\eqref{eqn:exp}.
\end{prop}

\pf  
By virtue of Theorem~\ref{T:1.1} and the fact that $r\to h_{ a, \gamma,  T}(t, r)$ is decreasing, the theorem holds for  $\delta_D(x) \ge at^{1/\alpha}/2$. 

We now fix $(t, x, y)\in (0, T]\times D\times D$ with $\delta_D(x)<at^{1/\alpha}/2$ and $|x-y|\ge  12at^{1/\alpha}{\bf 1}_{\b\in[0, 1]} +2 \cdot {\bf 1}_{\b\in(1, \infty)}+2(1+2at^{1/\alpha} )\cdot {\bf 1}_{\b= \infty}$, and we define $r_t:=at^{1/\alpha}$.
Let $U_1:=B(z_x, r_t)\cap D$, $U_3:=\{z\in D: |z-x|>|x-y|/2\}$, and $U_2:=D\backslash (U_1\cup U_3)$.
Then, $x\in U_1$ and $y\in U_3$. 
For $z\in U_2$, $|x-y|/2\le |x-y|-|x-z|\le |z-y|$.
Thus, by virtue of Theorem ~\ref{T:1.1}, we have
\begin{align*}
\sup_{s<t, z\in U_2} p(s,z,y)\le c_0\sup_{s<t,|z-y|>|x-y|/2} h_{ C_1, \gamma_1, T}(s, |z-y|) 
\le c_1\left( 1\vee (6a)^{-d-\alpha}\right)\, h_{ C_1, \gamma_1, T}(t, |x- y|/2).
\end{align*}
In fact, if $\b\in (1, \infty]$, we have $|z-y|\ge|x-y|/2>1$ and so $h_{C_1, \gamma_1, T}(s, |z-y|)$ is increasing in $s$. Also, if $\b \in [0, 1]$, we have $|z-y|\ge|x-y|/2\ge 6at^{1/\alpha}$ and $sr^{-\alpha-d}e^{-\gamma r^{\b}}$ is increasing in $s$. Thus, combining there observations with the fact $r\to h_{C_1, \gamma_1, T}(t, r)$ is decreasing, the second inequality above holds.

Moreover, from Lemma~\ref{L:2.3} and \eqref{e:L:3} in Theorem~\ref{L:2} we obtain
\begin{align}\label{tau}
\mathbb{P}_x\left(Y_{\tau_{U_1}}\in U_2\right)\le\mathbb{P}_x\left(Y_{\tau_{U_1}}\in D\right)\le c_2t^{-1}\mathbb{E}_x[\tau_{U_1}]\le c_3\frac{\delta_D(x)^{\alpha/2}}{\sqrt{t}} .
\end{align}
Hence, the first part of \eqref{eq:ub2} in Lemma~\ref{L:4.1} is bounded as follows:
\begin{align}\label{I}
\mathbb{P}_x\left(Y_{\tau_{U_1}}\in U_2\right) \Big(\sup_{s<t, z\in U_2} p(s, z, y)\Big)\le c_4 \frac{\delta_D(x)^{\alpha/2}}{\sqrt{t}}  h_{ C_1, \gamma_1, T}(t, |x-y|/2).                      
\end{align}

If $\b\in[0,\infty)$, since $|x-y|\ge 12at^{1/\alpha}$ we have for $u\in U_1$ and $z\in U_3$ that
\begin{align}\label{e:lbnew11}
|u-z|\ge|z-x|-|x-z_x|-|u-z_x|>  |x-y|/2-2at^{1/\alpha}\ge |x-y|/3.
\end{align}
Then, from \eqref{eqn:exp}--\eqref{e:conkappa1} and \eqref{e:L:3} we obtain
\begin{align}\label{II}
 \mathbb{E}_x[\tau_{U_1}]\Big(\sup_{u\in U_1, z\in U_3} J(u, z)\Big)
\le\, c_{5}\sqrt{t}{\delta_D(x)^{\alpha/2}}\frac{e^{-\gamma_1({|x-y|}/{3})^\b}}{|x-y|^{d+\alpha}} 
\le \,\, c_6\frac{\delta_D(x)^{\alpha/2}}{\sqrt{t}} h_{ \gamma_1, \gamma_1,T}(t, |x-y|/3).
\end{align}
If $\b=\infty$, since $|u-z| \ge |x-y|/2-2at^{1/\alpha}\ge1$ we have $J(u, z)=0$ on $U_1 \times U_3$. 
Hence, by applying \eqref{I} and \eqref{II} to \eqref{eq:ub2} for the case $\b\in[0, \infty)$ and applying \eqref{I} to \eqref{eq:ub2} for the case $\b=\infty$, we reach the conclusion.
\qed

For notational convenience, we denote by $X$ the process $Y$ in the case $\beta=0$, and we let $J^X(x, y):= \kappa(x,y){|x-y|^{-d-\alpha}}$ be its jumping kernel.
By Meyer's construction (e.g., see \cite[\S4.1]{CK2}), when $\beta \in (0, \infty]$ the process $Y$ can be constructed from $X$ by removing jumps of size greater than 1 with suitable rate.  
Let $p_D^X(t,x,y)$ be the transition density function of  $X$ on $D$. 
For $\beta \in (0, \infty]$, we define
\begin{align*}
\sJ(x):=\int_{\R^d} \kappa(x, y){|x-y|^{-(d+\alpha)}}\left(1- \psi_1(|x-y|)^{-1}\right)dy,
\end{align*}
where $\psi_1(|x-y|)$ is defined in \eqref{eqn:exp}.
Then, $\|\sJ\|_{\infty} \le c_1 \int_{|z| \ge 1 } |z|^{-(d+\alpha)}dz<\infty$.
By \cite[Lemma 3.6]{BBCK} we have 
\begin{align}
p_D(t,x,y)\,\le\, e^{T\|\sJ\|_{\infty}} p_D^X(t,x,y) \quad \mbox{ for any } (t,x,y) \in (0, T] \times D \times D. \label{e:u.0}
\end{align}
Thus, when $|x-y| <M$ for some $M>0$, it suffices to obtain the upper bound of $p_D^X(t, x, y)$, which is given next.  
\begin{prop}\label{P:X}
There exists a positive constant $c=c(R, \Lambda, \eta, \rho, T)$ such that for any $(t, x, y)\in (0, T]\times D\times D$ we have  
\[
p_D^X(t, x, y)\le c\Psi(t,x)\Psi(t,y)\left(t^{-d/\alpha}\wedge t|x-y|^{-\alpha-d}\right).
\]
\end{prop}
\pf
The semigroup property, Theorem~\ref{T:1.1} (for $\beta=0$), and Lemma~\ref{l:sdfnW} yield
\begin{align*}
p_D^X(t/2, x, y) \le \left(\sup_{z, w\in D} p_D^X(t/4, z, w)\right)\int_D p_D^X(t/4, x, z)dz
\le \,c_1t^{-d/\alpha} \P_x(\tau_D >t/4)  \le c_2 t^{-d/\alpha} \Psi(t,x).
\end{align*}
Thus, by Proposition~\ref{Pr:4.4} and Theorem~\ref{T:1.1} (for $\beta=0$), we obtain $p_D^X(t/2, x, y) \le  c_3 \Psi(t,x) p^X(t/2, x, y)$.
Combining this with Theorem~\ref{T:1.1} (for $\beta=0$), we conclude that 
\begin{align*}
&p_D^X(t, x, y)\,=\, \int_D p_D^X(t/2, x, z)\cdot p_D^X(t/2, z, y)dz\,\le \,c_3^2\Psi(t,x) \Psi(t,y)\int_{\R^d}  p^X(t/2, x, z)p^X(t/2, z, y)dz\\
&=\, c_3^2 \Psi(t,x) \Psi(t,y) p^X(t, x, y)\,\le \, c_{4} \Psi(t,x) \Psi(t,y) \left(t^{-d/\alpha}\wedge t|x- y|^{-\alpha-d}\right). 
\end{align*}
\qed

Combining Propositions~\ref{Pr:4.4} and~\ref{P:X}, we have the following proposition.

\begin{prop}\label{Pr:4.4n}
There exists a positive constant $c=c(\b, R, \Lambda, T, \eta, \rho)$ such that for every $(t, x, y)\in (0, T]\times D\times D$ we have 
\[
p_D(t, x, y)\le c \Psi(t,x)\cdot 
\begin{cases}
h_{ C_1\wedge\gamma_1,\gamma_1, T}(t, |x-y|/3)&\mbox{ if } \b\in[0,\infty),\\
h_{ C_1, \gamma_1, T}(t, |x-y|/2)&\mbox{ if } \b=\infty,
\end{cases}
\]
where $C_1$ is the constant in Theorem~\ref{T:1.1} and $\gamma_1$ is the constant in~\eqref{eqn:exp}.
\end{prop}

Next, we provide the upper bound estimates for $p_D(t, x, y)$ in the case $ \beta \in (0, \infty]$.

\medskip

\noindent {\bf Proof of Theorem~\ref{t:main}(1)}. 
By \eqref{e:u.0} and Proposition~\ref{P:X}, the theorem holds for $|x-y|\le  6(1 \vee C_1^{-1})$.
In fact, if $\beta = \infty$ and $6 <|x-y|\le  6(1 \vee C_1^{-1})$, then by \eqref{e:u.0} and Proposition~\ref{P:X} we have 
\begin{align}\label{e:4new1}
p_D(t, x, y)\le c_1 \Psi(t,x)\Psi(t,y) (t/T) \le c_1 \Psi(t,x)\Psi(t,y) (t/T)^{C_1|x-y|/6}.
\end{align}
The proofs for the other cases are obvious from \eqref{e:u.0} and Proposition~\ref{P:X}. 

Thus, by virtue of \eqref{e:u.0}, Proposition~\ref{P:X}, Proposition~\ref{Pr:4.4n}, and the symmetry of $p_D(t,x, y)$, we need to show the result only for the following case, which is assumed throughout the proof: $ |x-y|>  6(1 \vee C_1^{-1}) $ and $\delta_D(x)\vee \delta_D(y)<a_T t^{1/\alpha}$.

We define $r_t:=a_T t^{1/\alpha}$. 
For any $x$ with $\delta_D(x)<r_t$, let $z_x\in \partial D$ such that $\delta_D(x)=|z_x-x|$. Let $U_1:=B(z_x, r_t)\cap D$, $U_3:=\{z\in D: |z-x|>|x-y|/2\}$, and $U_2:=D\backslash (U_1\cup U_3)$.
Note that $x\in U_1$ and $y\in U_3$ and $|x-y|/2\le |z-y|$ for $z\in U_2$.
Thus, by Proposition~\ref{Pr:4.4n} we have   
\begin{align}\label{eq:p_D}
&\sup_{s<t, z\in U_2} p_D(s,z,y)\nn\\
\le&\sup_{s<t, z\in U_2}  c_2  \frac{\delta_D(y)^{\alpha/2}}{\sqrt{s}}\cdot 
\left(h_{ C_1\wedge \gamma_1, \gamma_1, T}(s, |z- y|/3)\cdot{\bf 1}_{\b\in[0, \infty)}
+h_{C_1, \gamma_1, T}(s, |z- y|/2)\cdot{\bf 1_{\b=\infty}}\right)\nn\\ 
\le & \,c_2 \,\delta_D(y)^{\alpha/2} \sup_{s<t, |x-y|/2\le |z-y|}\frac{1}{\sqrt{s} }\cdot
\left(h_{ C_1\wedge \gamma_1, \gamma_1, T}(s, |z- y|/3)\cdot{\bf 1}_{\b\in[0, \infty)}
+h_{C_1, \gamma_1, T}(s, |z- y|/2)\cdot{\bf 1}_{\b=\infty}\right)\nn\\ 
\le& \, c_3\,\frac{\delta_D(y)^{\alpha/2}}{\sqrt{t}}\cdot  \left(h_{ C_1\wedge \gamma_1, \gamma_1, T}(t, |x- y|/6)\cdot{\bf 1}_{\b\in[0, \infty)}
+h_{C_1, \gamma_1, T}(t, |x- y|/4)\cdot{\bf 1_{\b=\infty}}\right).
\end{align}
The last inequality is clear for $\beta \in [0, \infty)$ by definition of $h_{a, \gamma,T}$, and for $\beta=\infty$ we used the fact that $s \to s^{-1/2} (s/Tr)^{ar}$ is increasing if $ar \ge 1$. 
Hence, from \eqref{tau} and \eqref{eq:p_D} we obtain
\begin{align}\label{main:I}
&\mathbb{P}_x\left(Y_{\tau_{U_1}}\in U_2\right)\Big(\sup_{s<t, z\in U_2} p_D(s, z, y)\Big)\nn\\
&\le c_4 \frac{\delta_D(x)^{\alpha/2}}{\sqrt{t}} \frac{\delta_D(y)^{\alpha/2}}{\sqrt{t}} \cdot
\begin{cases}
h_{ C_1\wedge\gamma_1, \gamma_1, T}(t, |x-y|/6)&\mbox{ if } \b\in[0,\infty),\\
h_{ C_1, \gamma_1, T}(t, |x-y|/6)&\mbox{ if } \b=\infty.
\end{cases}
\end{align}

However, by Lemma~\ref{l:sdfnW} we have 
\begin{align}\label{main:II-1}
&\int_0^t \mathbb{P}_x(\tau_{U_1}>s) \mathbb{P}_y(\tau_D>t-s) ds \le \int_0^t \mathbb{P}_x(\tau_{D}>s) \mathbb{P}_y(\tau_D>t-s) ds \nn\\
&\qquad\le c_5 \, \delta_D(x)^{\alpha/2} \delta_D(y)^{\alpha/2}\int_0^t s^{-1/2}  (t-s)^{-1/2} ds \le c_6 \, t \,\frac{\delta_D(x)^{\alpha/2}}{\sqrt{t}} \frac{\delta_D(y)^{\alpha/2}}{\sqrt{t}}.
\end{align}
For $\b\in[0,\infty)$, we have $|u-z|\ge |x-y|/3$ for $(u,z) \in U_1 \times U_3$ as in \eqref{e:lbnew11}.
Thus, from \eqref{eqn:exp}--\eqref{e:conkappa1} and \eqref{main:II-1} we obtain 
\begin{align}\label{main:II}
 &\int_0^t \mathbb{P}_x\left(\tau_{U_1}>s\right) \mathbb{P}_y \left(\tau_D>t-s \right)ds \cdot \Big(\sup_{u\in U_1, z\in U_3} J(u, z)\Big)\nn\\
&\le c_{7}\frac{\delta_D(x)^{\alpha/2}}{\sqrt{t}}\frac{\delta_D(y)^{\alpha/2}}{\sqrt{t}}t\frac{e^{-\gamma_1({|x-y|}/{3})^\b}}{|x-y|^{d+\alpha} } \le c_{7}\frac{\delta_D(x)^{\alpha/2}}{\sqrt{t}}\frac{\delta_D(y)^{\alpha/2}}{\sqrt{t}}h_{ \gamma_1, \gamma_1,T}(t, |x-y|/3).
\end{align}
If $\b=\infty$, since $|u-z|> 1$, $J(u, z)=0$ on $ U_1 \times U_3$.
Therefore, by applying \eqref{main:I} and \eqref{main:II} in \eqref{eq:ub1} of Lemma~\ref{L:4.1} for $\b\in[0, \infty)$ and applying \eqref{main:I} for $\b=\infty$, we prove the theorem for $ |x-y|> 6(1 \vee C_1^{-1})$ and $\delta_D(x)\vee \delta_D(y)<a_T t^{1/\alpha}$.
\qed

\section{Lower bound estimates}

We proved the preliminary lower bound estimates in Section~3.
In this section, combining these results with the key estimate in \eqref{e:L:3-1}, we give the full lower bound estimate for $p_D(t, x, y)$ with the boundary decay terms.
We first introduce the next  lemma.

\begin{lemma}\label{L:5.1}
Suppose that $E_1, E_2, E$ are open subsets of $\bR^{d}$ with $E_1, E_2 \subset E$ and $\text{dist} (E_1, E_2)>0$.
If $x\in E_1$ and $y \in E_2$, then for all $t>0$ we have
\[
p_E(t,x,y)\ge t\,\ \bP_x(\tau_{E_1}>t)\,\ \bP_y(\tau_{E_2}>t) \inf_{(u, w)\in E_1\times E_2}J(u, w).
\]
\end{lemma}
\pf
See the proof of \cite[Lemma 3.3]{CKS7}. \qed

For the remainder of this section we assume that $\eta\in (\alpha/2, 1]$, $T>0$, and $D$ is a $C^{1,\eta}$ open set with characteristics $(R, \Lambda)$.
Without loss of the generality, we assume that $\Lambda>4$ and $R <10^{-1}$.
We let 
\begin{equation*}
\wh a_T=a_{T,R}:=2^{-5} {R}{T^{-1/\alpha}} <2^{-5}10^{-1} T^{-1/\alpha},
\end{equation*}
and for $x\in D$ we use $z_x$ to denote a point on $\partial D$ such that  $|z_x-x|=\delta_D(x)$.

The next two lemmas are crucial to obtain the lower bound on the survival probability where $x$ is near the boundary of $D$.

\begin{lemma}\label{L:5.2}
For any $a \le \wh a_T$, there exists a positive constant $c=c (a, \b, R, \Lambda, T, \eta, \rho)$ such that for every $t<T$ and  $x\in D$ with $ \delta_{D} (x)<at^{1/\alpha}$ we have
\[
\bP_x(\tau_{B(z_x, 10at^{1/\alpha})\cap D}>t/3)\ge c \,\ \frac{\delta_{D} (x)^{\alpha/2}}{\sqrt{t}}.
\]
\end{lemma}

\pf
Without loss of generality, we assume that $z_x=0$.
Consider a coordinate system $CS:=CS_0$ such that $B(0, R)\cap D=\{y=(\wt{y}, y_d)\in B(0, R) \text{ in } CS: y_d>\phi(\wt{y})\}$, where $\phi$ is a $C^{1,\eta}$ function such that $\phi(0)=0$, $ \nabla\phi (0)=(0, \dots, 0)$, $\| \nabla\phi \|_\infty \leq \Lambda$, and $| \nabla \phi (\wt y)-\nabla \phi (\wt w)| \leq \Lambda |\wt y-\wt w|^{\eta}$.

Let $\psi(\wt{y})=2\Lambda|\wt{y}|$ and $V:=\{y=(\wt{y}, y_d)\in B(0, R) \text{ in } CS: y_d>\psi(\wt{y})\}$.
Then, since $\psi(\wt{y}) \ge 2\Lambda|\wt{y}|^{\eta+1}$, the mean value theorem yields $\{y=(\wt{y}, y_d)\in B(0, R) \text{ in } CS: y_d>\psi(\wt{y})\} \subset B(0, R)\cap D$.

Let $U_1:=B(0, 2at^{1/\alpha})\cap D$, $U_2:=B(0, 10at^{1/\alpha})\cap D$, and   
\[
W:=\{y=(\wt{y}, y_d)\in B(0,8 at^{1/\alpha})\backslash B(0, 2at^{1/\alpha}) \text{ in } CS: y_d>\psi(\wt{y})\}.
\]
Since $ \Lambda |\wt{w}|=\psi(\wt{w})/2<w_d /2$ for $w\in W$, we have
\begin{align}\label{e:jk1}
 \delta_D (w)> (1+\Lambda)^{-1}(w_d-\phi(\wt{w}))> (1+\Lambda)^{-1}(w_d-\Lambda |\wt{w}|)> 2^{-1}(1+\Lambda)^{-1} w_d.
\end{align}
Moreover, since $|\wt w| \le (2\Lambda)^{-1}   |w| \le \Lambda^{-1}4a t^{1/\alpha}\le a t^{1/\alpha} $ for $w\in W$, 
we have
\begin{align}\label{e:jk2} 
w_d^2> (2at^{1/\alpha})^2- |\wt w|^2 >2at^{1/\alpha} (2at^{1/\alpha}- |\wt w|) 
\ge (at^{1/\alpha})^2
\text{ for  } w\in W.
\end{align}
Combining \eqref{e:jk1} and \eqref{e:jk2}, we obtain $\delta_D (w)>  2^{-1}(1+\Lambda)^{-1}at^{1/\alpha}$.
Thus, $B(w,r_1at^{1/\alpha})\subset U_2$ for $w\in W$, where $r_1:=2^{-1}(1+\Lambda)^{-1}$.
Hence, by virtue of the strong Markov property, Lemma~\ref{L:3.3}, and \eqref{e:L:3-1}, we have
\begin{align*}
&\bP_x(\tau_{U_2}>t/3)\,\geq\, \bP_x(\tau_{U_2}>t/3, Y_{\tau_{U_1}} \in W)\,=\, \bE_x[\bP_{Y_{\tau_{U_1}}}(\tau_{U_2}>t/3): Y_{\tau_{U_1}}\in W]\\
&\geq \bE_x[\bP_{Y_{\tau_{U_1}}}(\tau_{B(Y_{\tau_{U_1}}, r_1at^{1/\alpha})}>t/3): Y_{\tau_{U_1}}\in{W}]\\
&= \bP_0(\tau_{B(0, r_1at^{1/\alpha})}>t/3)\bP_x(Y_{\tau_{U_1}}\in{W}) \geq c_1\,\bP_x(Y_{\tau_{U_1}}\in{W})
\geq c_2 \frac{\delta_D(x)^{\alpha/2}}{\sqrt{t}}.
\end{align*}
\qed

We introduce the following definition for the subsequent lemma.
\begin{defn}\label{def:UB}
Let $0<\kappa\leq 1/2$.
We say that an open set $D$ is $\kappa$-fat if there is $R_1>0$ such that for all $x\in \overline{D}$ and all $r\in (0,R_1]$ there is a ball $B(A_r(x), \kappa r)\subset D\cap B(x,r)$.
The pair $(R_1, \kappa)$ are called the characteristics of the $\kappa$-fat open set $D$.
\end{defn}
It is clear that a $C^{1, \eta}$ open set $D$ with characteristics $(R, \Lambda)$ is always a $\kappa$-fat set whose characteristics $(R_1, \kappa)$ depend only on $R$, $\Lambda$, and $d$.
Hereinafter, without loss of generality, we assume that $R \le R_1$ (by choosing $R$ smaller if necessary) and that $A_r(x)$ is always the point $A_r(x)\in D$ in Definition~\ref{def:UB} for $D$.
Recall that $\Psi$ is defined in~\eqref{e:dax}.

\begin{lemma}\label{L:5.3}
Let $T$ be a positive constant.
For any $\b \in [0, \infty]$, there exists a positive constant $c=c(\b, R, \Lambda, T, \eta, \rho)>0$ such that, for every $t < T$ and $x\in D$, we can find $x_1$ with $\delta_D(x_1)\ge  2^{-1}\kappa \wh a_Tt^{1/\alpha}$  and $|x_1-x|  \le  6\wh a_Tt^{1/\alpha}$ such that 
\[
\int_{B(x_1, (\kappa/4) 
 \wh a_Tt^{1/\alpha})} p_D(t/3, x, z)dz \ge c \Psi(t,x). 
\]
\end{lemma}

\pf
For $\delta_D(x)< 2^{-1}\kappa \wh a_Tt^{1/\alpha}$, let $x_1=A_{6\wh a_Tt^{1/\alpha}}(z_x)$. 
Let $B_{x_1}:=B(x_1, (\kappa/4)  \wh a_Tt^{1/\alpha})$ and $B_{z_x}:=B(z_x, 5\kappa \wh a_Tt^{1/\alpha})\cap D$ so that $B_{x_1}\cap B_{z_x} =\emptyset$.
By Lemmas~\ref{L:5.1}, \ref{L:5.2}, and \ref{L:3.3}, \begin{eqnarray*}
\int_{B_{x_1}} p_D(t/3,x,z)dz 
&\ge& \frac{t}{3} \int_{B_{x_1}} \bP_x(\tau_{B_{z_x}}>t/3)\,\,\bP_z(\tau_{B_{x_1}}>t/3)  \cdot \inf_{(u, w)\in B_{z_x}\times B_{x_1}}J(u, w) dz\\
&\ge& \frac{t}{3}\,\ \bP_x(\tau_{B_{z_x}}>t/3)\cdot c_1 \int_{B_{x_1}} dz\cdot c_2 \frac{1}{(t^{1/\alpha})^{d+\alpha}}\\
&=&c_3 \bP_x(\tau_{B_{z_x}}>t/3)\,\ge\, c_4\,\ \frac{\delta_D(x)^{\alpha/2}}{\sqrt{t}}.
\end{eqnarray*}
For $\delta_D(x)\ge 2^{-1}\kappa \wh a_Tt^{1/\alpha}$, let $x_1=x$ and $B_{x_1}:=B(x_1,  (\kappa/4)  \wh a_Tt^{1/\alpha})$.
By Lemma~\ref{L:3.3}, there exists a constant $c_5=c_5(\alpha,\b, R, T, d,L_3)>0$ such that
\begin{eqnarray*}
\int_{B_{x_1}} p_D(t/3,x,z)dz 
\ge \int_{B_{x_1}} p_{B_{x_1}}(t/3, x, z) dz= \bP_x(\tau_{B_{x_1}}>t/3)> c_5 .
\end{eqnarray*}
This proves the lemma.
\qed

We are now ready to give the proof of the lower bound estimates for $p_D(t, x, y)$.
Recall our assumption that $\eta\in (\alpha/2, 1]$ and $D$ is a $C^{1,\eta}$ open set. 
For the cases $\b\in(1,\infty)$ with $|x-y| \ge 1$ and $\b=\infty$ with $|x-y| > 4/5$, we assume in addition that the path distance in each connected component of $D$ is comparable to the Euclidean distance with characteristic $\lambda_1$. 
Note that combining this assumption with $C^{1,\eta}$ assumption entails that $D$ satisfies the assumption made before Proposition~\ref{step5}.  

\medskip

\noindent {\bf Proof of Theorem \ref{t:main}(2) and \ref{t:main}(3)}. 
By Lemma~\ref{L:5.3}, for any $x, y\in D$, there exists $x_1, y_1\in D$ such that $\delta_D(x_1) \wedge \delta_D(y_1) \ge 2^{-1}\kappa \wh a_Tt^{1/\alpha}$ and $|x_1-x|  \vee |y_1-y| \le  6\wh a_Tt^{1/\alpha}$, and there exists a constant $c_1=c_1(\eta, \rho, \b, R, \Lambda, T)>0$ independent of $x,y$ such that
\begin{align}\label{eq:m11}
\int_{B_{x_1}} p_D(t/3, x, z)dz\int_{B_{y_1}} p_D(t/3, y, z)dz \ge c_1 \Psi(t,x) \Psi(t,y),
\end{align}
where $B_{x_1}:=B(x_1, (\kappa/4)  \wh a_Tt^{1/\alpha})$ and $B_{y_1}:=B(y_1, (\kappa/4)  \wh a_Tt^{1/\alpha})$.
Thus, by the semigroup property we have
\begin{align}\label{eq:m1}
p_D(t, x, y)
&=\int_D \int_D p_D(t/3, x, u) p_D(t/3, u, w) p_D(t/3,w,y)du dw\nonumber\\ 
&\ge \int_{B_{x_1}} p_D(t/3,x,u)du \int_{B_{y_1}} p_D(t/3, y, w) dw  \left(\inf_{(u,w)\in B_{x_1}\times B_{y_1}} p_D(t/3, u, w) \right) \nonumber\\
&\ge c_1\Psi(t,x) \Psi(t,y) \inf_{(u,w)\in B_{x_1}\times B_{y_1}}  p_D(t/3, u, w) .
\end{align}

We now carefully calculate the lower bounds of $p_D(t/3,u,w)$ on $B_{x_1} \times B_{y_1}$.
Since $|x-x_1|\vee |y-y_1|\le 6\wh a_Tt^{1/\alpha}$, for $u\in B_{x_1}$ and $w\in B_{y_1}$ we have
\begin{align}\label{e:dfdf5}
&|x-y|-20^{-1} \le |x-y|-(12+ (\kappa/2))\wh a_Tt^{1/\alpha} \nn\\
&\le |u-w| \le |x-y|+(12+ (\kappa/2))\wh a_Tt^{1/\alpha} \le  |x-y|+ 20^{-1}
\end{align}
and 
$\delta_D(u)\wedge \delta_D(w)\ge (\kappa/4) \wh a_Tt^{1/\alpha}$. 

If $\b\in[0,1]$, then by considering the cases $|x-y| <15 \wh a_Tt^{1/\alpha}$ and $|x-y| >15 \wh a_Tt^{1/\alpha}$ separately using Proposition~\ref{step4} and \eqref{e:dfdf5} we obtain
\begin{align*}
p_D(t/3, u, w)\ge c_2 \left(t^{-d/\alpha} \wedge t|u-w|^{-d-\alpha} e^{-\gamma_2|u-w|^\b} \right)
\ge c_3 \left(t^{-d/\alpha} \wedge t|x-y|^{-d-\alpha} e^{-\gamma_2|x-y|^\b} \right) .
\end{align*}

If $\b\in(1, \infty]$ and $|x-y| \le 4/5$, then \eqref{e:dfdf5} yields $|u-w|\le|x-y|+20^{-1} <1$. 
Thus, by considering the cases $|x-y| <15 \wh a_Tt^{1/\alpha}$ and $|x-y| >15 \wh a_Tt^{1/\alpha}$ separately using Proposition~\ref{step4} and \eqref{e:dfdf5}, we have $p_D(t/3, u, w)\ge c_4 \left(t^{-d/\alpha} \wedge (t|x-y|^{-d-\alpha}) \right)$.

If $\b\in(1, \infty]$ and $4/5\le |x-y|$, then \eqref{e:dfdf5} yields $|u-w| \asymp |x-y|$.

We now consider $p_D(t/3, u, w)$ in each of the remaining cases. 
\begin{itemize}
\item[(1)]
If $\b\in(1, \infty)$ and $4/5 \le |x-y|<2$, then $|u-w|\asymp 1$.
Thus, by Proposition~\ref{step4}, we have $p_D(t/3, u, w)\ge c_5t$.

\item[(2)]
If $\b=\infty$ and $4/5\le |x-y|<2$, then by Propositions~\ref{step4} and~\ref{step5} we have
\begin{align*}
p_D(t/3, u,w)\ge c_6 \frac{4t }{5T|x-y|}\ge c_{6}\left(\frac{4t}{5T|x-y|}\right)^{5|x-y|/4}.
\end{align*}

\item[(3)]
If $\b\in(1, \infty)$ and $2\le |x-y|$, then $1 <|u-w|$ and from Proposition~\ref{step6} and \eqref{e:dfdf5} we obtain
\begin{align*}
&p_D(t/3, u, w)\ge c_7  t\exp\left\{-c_{8}\left(|u-w|\left(\log\frac{T|u-w|}{t}\right)^{\frac{\b-1}{\b}}\wedge |u-w|^{\b}\right)\right\}\\
&\ge c_7  t\exp\left\{-c_{8}\left((5|x-y|/4)\left(\log\left(\frac{T(|x-y|+20^{-1})}{t}\right)\right)^{\frac{\b-1}{\b}}\wedge (5|x-y|/4)^{\b}\right)\right\}\\
&\ge c_7  t\exp\left\{-c_{9}\left(|x-y|\left(\log\frac{T|x-y|}{t}\right)^{\frac{\b-1}{\b}}\wedge |x-y|^{\b}\right)\right\}.
\end{align*}
The last inequality comes from the inequality $\log r \le \log(r+b)\le 2\log r$ for $r\ge 2\vee b$.  

\item[(4)]
If $\b=\infty$ and $2\le |x-y|$, then $1<|u-w|$ and from Proposition~\ref{step5} and \eqref{e:dfdf5} we have 
\begin{align*}
&p_D(t/3, u, w)\ge c_{11}\left(\frac{t}{T|u-w|}\right)^{c_{10} |u-w|}
\ge c_{11}\left(\frac{t}{T(|x-y|+20^{-1})}\right)^{c_{12}|x-y|}\\
&\ge c_{11}\left(\frac{t}{T|x-y|}\right)^{2c_{12}|x-y|}\ge c_{11}\left(\frac{4t}{5T|x-y|}\right)^{2c_{12}|x-y|}.
\end{align*}
The second last inequality holds by virtue of the inequality $r^2\ge r+b$ for $r\ge 2\vee b$.
\end{itemize}

Hence, combining \eqref{eq:m1} with the above observations on the lower bound of $p_D(t/3, u, w)$, we have proved Theorem~\ref{t:main}(2) and~\ref{t:main}(3).
\qed

\medskip

\noindent {\bf Proof of Theorem \ref{t:main}(4)}. 
Let $D(x)$ and $D(y)$ be connected components containing $x$ and $y$, respectively.
By definition of a $C^{1, \eta}$ open set, the distance between $x$ and $y$ is at least $R$.
Using Lemma~\ref{L:5.3}, we find that $x_1\in D(x)$ and $y_1\in D(y)$.
We then define $B_{x_1}$ and $B_{y_1}$ in the same way as when beginning the proof of Theorem~\ref{t:main}(2) and~\ref{t:main}(3) so that \eqref{eq:m11} holds and for any $u\in B_{x_1}$ and $w\in B_{y_1}$ we have $3R/4 \le 3|x-y|/4\le  |u-w|\le 5|x-y|/4.$
By Proposition~\ref{step4}, for every $u\in B_{x_1}$ and $w\in B_{y_1}$ we have
\[
p_D(t/3, u, w)\ge  c_{1} \frac{t }{ |u-w|^{d+\alpha}}e^{-\gamma_2 |u-w|^\b}\ge c_2\frac{t}{ |x-y|^{d+\alpha}}e^{-\gamma_2(5|x-y|/4)^\b}.
\]
Therefore, 
\begin{align*}
p_D(t, x, y)
&\ge\int_{B_{x_1}} \int_{B_{y_1}} p_D(t/3, x, w) p_D(t/3, u, w) p_D(t/3,w,y)dw dv\nonumber\\ 
&\ge \int_{B_{x_1}} p_D(t/3,x,u)du \int_{B_{y_1}} p_D(t/3, y, w) dw \cdot \inf_{(u, w)\in B_{x_1} \times B_{y_1}} p_D(t/3, u, w) \nonumber\\
&\ge c_3 \Psi(t,x) \Psi(t,y) \cdot \frac{t}{ |x-y|^{d+\alpha}}e^{-\gamma_2(5|x-y|/4)^\b}.
\end{align*}
\qed

\medskip

\noindent {\bf Proof of Theorem \ref{t:main}(5)}. 
Note that, since $D$ is bounded and connected, the estimate for $p_D(t, x, y)$ at small time is the same as that obtained for a symmetric stable process in \cite{CKS}.
Thus, the remainder of the proof of Theorem~\ref{t:main}(5) using the estimate for $p_D(t, x, y)$ at small time is routine (see \cite{CKS}) and we omit it here.

\begin{small}

\end{small}

\vskip 0.3truein

{\bf Kyung-Youn Kim}

Department of Mathematical Sciences,
Seoul National University,
Building 27, 1 Gwanak-ro, Gwanak-gu,
Seoul 151-747, Republic of Korea

E-mail: \texttt{eunicekim@snu.ac.kr}

\bigskip

{\bf Panki Kim}

Department of Mathematical Sciences and Research Institute of Mathematics,
Seoul National University,
Building 27, 1 Gwanak-ro, Gwanak-gu,
Seoul 151-747, Republic of Korea

E-mail: \texttt{pkim@snu.ac.kr}

\end{document}